\documentclass[12pt]{amsart}

\usepackage{mathrsfs}

\usepackage{amsfonts, amssymb,amsmath, amsthm, amsxtra, latexsym, amscd}
\usepackage[latin1]{inputenc} 

\usepackage[english]{babel}

\usepackage{amsmath}

\usepackage{amsthm,amssymb,latexsym}

\usepackage{t1enc}

\usepackage{epic}

\usepackage{eepic}

\usepackage{xypic}

\usepackage{amsfonts, amssymb,amsmath, amsthm, amsxtra, latexsym, amscd}

\usepackage{graphics}

\title[Branching laws for minimal holomorphic representations]{Branching laws for
minimal holomorphic representations}
\author{Henrik Seppänen}
\address{Department of Mathematics, Chalmers University of Technology 
and G\"oteborg University, G\"oteborg, Sweden}
\email{henriks@math.chalmers.se}

\newtheorem{prop}{Proposition}

\newtheorem{lemma}[prop]{Lemma}

\newtheorem{cor}[prop]{Corollary}

\newtheorem{thm}[prop]{Theorem}

\theoremstyle{definition}

\newtheorem{defin}[prop]{Definition}

\theoremstyle{remark}

\begin{document}





\maketitle

\begin{abstract}
In this paper we study the branching law for the restriction from $SU(n,m)$ to $SO(n,m)$ of the minimal representation in the
analytic continuation of the scalar holomorphic discrete series. We identify the the group decomposition
with the spectral decomposition of the action of the Casimir operator on the subspace of $S(O(n) \times O(m))$-invariants.
The Plancherel measure of the decomposition defines an $L^2$-space of functions, for which certain continuous dual Hahn polynomials
furnish an orthonormal basis. It turns out that the measure has point masses precisely when $n-m>2$. Under these conditions
we construct an irreducible representation of $SO(n,m)$, identify it with a parabolically induced representation, and construct
a unitary embedding into the representation space for the minimal representation of $SU(n,m)$.
\end{abstract}

\section{Introduction}
One of the most important problems in harmonic analysis and in representation theory is that of decomposing group representations into
irreducible ones. When the given representation  
arises as the restriction of an irreducible representation of a bigger group, the decomposition is referred to as a \emph{branching law}.
One of the most famous examples of this is the Clebsch-Gordan decomposition for the restriction of the tensor product of two irreducible
$SU(2)$-representations (which is a representation of $SU(2) \times SU(2))$ to the diagonal subgroup.
For an introduction to the general theory for compact connected Lie groups, we refer to \cite{knapp-beyond}.

Since the work by Howe (\cite{howedual}) and by Kashiwara-Vergne (\cite{k-w}), the study of branching rules for singular and minimal 
representations on spaces of holomorphic functions on bounded symmetric domains has been an active area of research.
In \cite{jak-v-restr}, Jakobsen and Vergne studied the restriction to the diagonal subgroup of two holomorphic representations.
More recently, Peng and Zhang (\cite{peng-zhang}) studied the corresponding decomposition for the tensor product 
of arbitrary (projective) representations in the analytic continuation of the scalar holomorphic discrete series. Zhang also studied
the restriction to the diagonal of a minimal representation in this family tensored with its own anti-linear dual (\cite{zhtp}).

The restriction of the representations given by the analytic continuation of the scalar holomorphic discrete series to symmetric
subgroups (fixed point groups for involutions) has been studied recently by Neretin (\cite{ner}, \cite{nerbeta}), 
Davidson, Ólafsson, and Zhang (\cite{doz}), 
Zhang (\cite{zhtams}, \cite{zhSB}), van Dijk and Pevzner (\cite{dijkpevz}) and by the author (\cite{papper1}).

All the above mentioned decompositions have the common feature that they are multiplicity free. This general result
follows from a recent theorem by Kobayashi (\cite{kob-prop-vb}), where some geometric conditions are given for the action of a Lie group
as isometric automorphisms of a Hermitian holomorphic vector bundle over a connected complex manifold to 
guarantee the multiplicity-freeness in the decomposition 
of any Hilbert space of holomorphic sections of the bundle. The action of a symmetric subgroup on the trivial line bundle over
a bounded symmetric domain
then satisfies these conditions (cf. \cite{kob-vis-symspace}).

In this paper we study the branching rule for the restriction from $G:=SU(n,m)$ to $H:=SO(n,m)$ of the minimal representation in the 
analytic continuation of the scalar holomorphic discrete series. We consider the the subspace of $L:=S(O(n) \times O(m))$-invariants
and study the spectral decomposition for the action of Casimir element of the Lie algebra of  $H$. The diagonalisation
gives a unitary isomorphism between the subspace of $L$-invariants and an $L^2$-space with a Hilbert basis
given by certain continuous dual Hahn polynomials. The main theorem is Theorem \ref{irruppdelning}, 
where the decomposition on the group level is identified with
this spectral decomposition. The Plancherel measure turns out to have point masses precisely when $n-m>2$. The second half of the paper is 
devoted to the realisation of the representation associated with one of these points and the unitary embedding into
the representation space for the minimal representation. The main theorem of the second half is Theorem \ref{embedding}.

The paper is organised as follows. In Section 2 we begin with some preliminaries on the structure of the Lie algebra $\mathfrak{g}$, the
group action, and the minimal representation. In Section 3 we construct an orthonormal basis for the subspace of $L$-invariants.
In Section 4 we compute the action of the Casimir elements on the $L$-invariants and find its diagonalisation. We also state the
branching theorem. In Section 5 we construct an irreducible representation of the group $H$ (for $n-m>2$, i.e.,  when
point masses occur in the Plancherel measure), identify it with a parabolically induced
representation, and finally we construct a unitary embedding that realises one of the discrete points in the spectrum.

{\bf Acknowledgement:} The author would like to thank his advisor Professor Genkai Zhang for support and for many valuable suggestions during the
preparation of this paper.

\section{Preliminaries}

Let $\mathscr{D}$ be the bounded symmetric domain of type $I_{mn} (n \geq m)$ , i.e., 
\begin{eqnarray}
\mathscr{D}:=\left\{z \in M_{nm}(\mathbb{C}) | I_n-zz^* > 0\right\}.
\end{eqnarray}
Here $M_{nm}(\mathbb{C})$ denotes the complex vector space of $n \times m$ matrices.
We let $G$ be the group $SU(n,m)$, i.e., the group of all complex $(n+m) \times (n+m)$ matrices of determinant 
one preserving the sesquilinear form $\langle \cdot, \cdot \rangle_{n,m}$ on $\mathbb{C}^{n+m}$ given by
\begin{eqnarray}
\langle u, v \rangle_{n,m}=u_1\overline{v}_1+ \cdots +u_n\overline{v}_n-
u_{n+1}\overline{v}_{n+1}-\cdots-u_{n+m}\overline{v}_{n+m}.
\end{eqnarray}
The group $G$ acts holomorphically on $\mathscr{D}$ by 
\begin{eqnarray}
g(z)=(Az+B)(Cz+D)^{-1}, \label{fraktionell G}
\end{eqnarray}
if $g=\left( \begin{array}{cc}
A & B \\
C & D \\
\end{array}\right)$ is a block matrix determined by the size of $A$ being $n \times n$.
The isotropy group of the origin is
\begin{eqnarray*}
K&:=&S(U(n) \times U(m))\\
&=&\left\{ \left( \begin{array}{cc}
A & 0\\
0 & D \\
\end{array}\right) | A \in U(n), D \in U(m), \det(A)\det(D)=1\right\}, \nonumber
\end{eqnarray*}
and hence
\begin{eqnarray}
\mathscr{D} \cong G/K.
\end{eqnarray}

\subsection{Harish-Chandra decomposition}
Let $\theta$ denote the Cartan involution $g \mapsto (g^*)^{-1}$ on $G$. We use the same letter to denote its differential 
$\theta: \mathfrak{g} \rightarrow \mathfrak{g}$ at the identity. Here, we have identified $T_e(G)$ with $\mathfrak{g}$.
Let 
\begin{eqnarray}
\mathfrak{g}=\mathfrak{k} \oplus \mathfrak{p}
\end{eqnarray}
be the decomposition into the $\pm 1$ eigenspaces of $\theta$ respectively. 
In terms of matrices, 
\begin{eqnarray}
\mathfrak{k}&=&\left\{ \left( \begin{array}{cc}
A & 0\\
0 & D \\
\end{array}\right)|A^*=-A, D^*=-D, \mbox{tr}(A)+\mbox{tr}(D)=0\right\},\\
\mathfrak{p}&=&\left\{ \left( \begin{array}{cc}
0 & B\\
B^* & 0 \\
\end{array}\right)\right\},
\end{eqnarray}
where the size $A$ is $n \times n$.

The Lie algebra $\mathfrak{g}$ has a compact Cartan subalgebra $\mathfrak{t} \subset \mathfrak{k}$, where
\begin{eqnarray}
&&\mathfrak{t}
=\left\{ \left( \begin{array}{cccccc}
is_1 & 0 & \cdots & 0 & \cdots & 0\\
\vdots & \ddots & \vdots & \vdots & \vdots & \vdots \\
0 & \cdots & is_n & 0 & \cdots & 0\\
0 & \cdots & 0 & it_1 & \cdots & 0\\
\vdots & \vdots & \vdots & \vdots & \ddots & \vdots\\
0 & \cdots & 0 & 0 & \cdots & it_m\\
\end{array}\right)| 
\begin{array}{c}
s_i, t_j \in \mathbb{R}\\
\sum_is_i+\sum_jt_j=0
\end{array}
\right\}.  \label{cartan}
\end{eqnarray}
Its complexification, $\mathfrak{t}^{\mathbb{C}}$ (the set of complex diagonal traceless matrices),
is a Cartan subalgebra of the complexification $\mathfrak{g}^{\mathbb{C}}=\mathfrak{sl}(n+m, \mathbb{C})$, where
\begin{eqnarray}
\mathfrak{g}^{\mathbb{C}}=\mathfrak{k}^{\mathbb{C}} \oplus \mathfrak{p}^{\mathbb{C}}.
\end{eqnarray}
We let $E_{ij}$ denote the matrix with $1$ at the entry corresponding to the $i$th row and the $j$th column and zeros elsewhere.
By $E_{ij}^*$ we mean the dual linear functional, i.e., $E_{ij}^*(z)=z_{ij}$ for $z \in M_{nm}(\mathbb{C})$.
Moreover, we define an ordered basis $\{F_j\}$ for $\mathfrak{t}^{\mathbb{C}}$ by
\begin{eqnarray}
&&F_j:=E_{jj}^*-E_{j+1\,j+1}^*, \,j=1, \ldots, n+m-1, \label{ordn}\\
&&F_1 \leq \cdots \leq F_{n+m-1}. \nonumber
\end{eqnarray}

The root system, $\Delta(\mathfrak{g}^{\mathbb{C}}, \mathfrak{t}^{\mathbb{C}})$ is given by
\begin{eqnarray}
\Delta(\mathfrak{g}^{\mathbb{C}}, \mathfrak{t}^{\mathbb{C}})=\{E^*_{ii}-E^*_{jj}|1 \leq i,j \leq n+m, i \neq j\}.
\end{eqnarray}
We denote the root $E^*_{ii}-E^*_{jj}$ by $\alpha_{ij}$.
We define a system of positive roots, $\Delta^+$, by the ordering \eqref{ordn}. Then 
\begin{eqnarray}
\Delta^+=\{\alpha_{ij}|j>i\},
\end{eqnarray}
and we let $\Delta^-$ denote the complement so that $\Delta=\Delta^+ \cup \Delta^-$.
For a root, $\alpha$, we let $\mathfrak{g}^{\alpha}$ stand for the corresponding root space. Then 
$\mathfrak{g}^{\alpha_{ij}}=\mathbb{C}E_{ij}$.
For a root space, $\mathfrak{g}^{\alpha}$, we either have $\mathfrak{g}^{\alpha} \subset \mathfrak{k}^{\mathbb{C}}$ or
$\mathfrak{g}^{\alpha} \subset \mathfrak{p}^{\mathbb{C}}$. In the first case, we call the corresponding root compact, and
in the second case we call it non-compact. We denote the sets of compact and non-compact roots by
$\Delta_{\mathfrak{k}}$ and $\Delta_{\mathfrak{p}}$ respectively.
Finally, we let $\Delta_{\mathfrak{p}}^+$ and $\Delta_{\mathfrak{p}}^-$ denote the set of non-compact positive roots
and the set of non-compact negative roots respectively. We set
\begin{eqnarray}
\mathfrak{p}^+&=&\sum_{\alpha \in \Delta_{\mathfrak{p}}^+}\mathfrak{g}^{\alpha},\\
\mathfrak{p}^-&=&\sum_{\alpha \in \Delta_{\mathfrak{p}}^-}\mathfrak{g}^{\alpha}.
\end{eqnarray}
These subspace are abelian Lie subalgebras of $\mathfrak{p}^{\mathbb{C}}$. Moreover, the relations
\begin{eqnarray}
[\mathfrak{k}^{\mathbb{C}}, \mathfrak{p}^+]\subseteq \mathfrak{p}^+, 
[\mathfrak{k}^{\mathbb{C}}, \mathfrak{p}^-]\subseteq\mathfrak{p}^-,
[\mathfrak{p}^+, \mathfrak{p}^-]\subseteq \mathfrak{k}^{\mathbb{C}}
\end{eqnarray}
hold. We let $K^{\mathbb{C}}, P^+$, and $P^-$ denote the connected Lie subgroups of the complexification of 
$G$, $G^{\mathbb{C}}$, with Lie algebras
$\mathfrak{k}^{\mathbb{C}}, \mathfrak{p}^+$, and $\mathfrak{p}^-$ respectively.
The exponential mapping $\exp: \mathfrak{p}^{\pm} \rightarrow P^{\pm}$ is a diffeomorphic isomorphism of abelian groups. 
As subspaces of the Lie algebra $\mathfrak{g}^{\mathbb{C}}=\mathfrak{sl}(n+m)$ we have the matrix realisations
\begin{eqnarray}
\mathfrak{p}^+&=&\left\{\left( \begin{array}{cc}
0 & z\\
0 & 0\\
\end{array}\right)|z \in M_{nm}(\mathbb{C})\right\},\\
\mathfrak{p}^-&=&\left\{\left( \begin{array}{cc}
0 & 0\\
z & 0\\
\end{array}\right)|z \in M_{mn}(\mathbb{C})\right\}.
\end{eqnarray}
The Lie algebra $\mathfrak{g}^{\mathbb{C}}$ can be decomposed as 
\begin{eqnarray}
\mathfrak{g}^{\mathbb{C}}=\mathfrak{p}^+ \oplus \mathfrak{k}^{\mathbb{C}} \oplus \mathfrak{p}^-.
\end{eqnarray}
On a group level, the multiplication map
\begin{eqnarray}
P^+ \times K^{\mathbb{C}} \times P^- \rightarrow G^{\mathbb{C}},
(p,k,q) \mapsto pkq
\end{eqnarray}
is injective, holomorphic and regular with open image containing $GP^+$. In fact, identifying
the domain $\mathscr{D}$ with the subset\\
$\left\{ \left( \begin{array}{cc}
0 & z\\
0 & 0\\
\end{array}\right)| z \in \mathscr{D}\right\} \subset \mathfrak{p}^+$ and letting 
\begin{eqnarray*}
\Omega:=\exp \mathscr{D}=
\left\{ \left( \begin{array}{cc}
I_n & z\\
0 & I_m\\
\end{array}\right)| z \in \mathscr{D}\right\},
\end{eqnarray*}
there is an inclusion
\begin{eqnarray}
GP^+ \subset \Omega K^{\mathbb{C}}P^-.
\end{eqnarray}
For $g \in G$, we let $(g)_+, (g)_0$, and $(g)_-$ denote its $P^+, K^{\mathbb{C}}$, and $P^-$ factors respectively.
The action of $g$ on $\mathscr{D}$ defined by
\begin{eqnarray}
g(z)=\log ((g \exp z)_+)
\end{eqnarray}
then coincides with the action \eqref{fraktionell G}. In fact, for $g=\left( \begin{array}{cc}
A & B \\
C & D \\
\end{array}\right)$, the Harish-Chandra factorisation is given by
\begin{eqnarray}
&&\left( \begin{array}{cc}
A & B \\
C & D \\
\end{array}\right) \label{konkret-HC}\\
&&=\left( \begin{array}{cc}
I_n & BD^{-1} \\
0 & I_m \\
\end{array}\right)
\left( \begin{array}{cc}
A-BD^{-1}C & 0 \\
0 & D \\
\end{array}\right)
\left( \begin{array}{cc}
I_n & 0 \\
D^{-1}C & I_m \\
\end{array}\right). \nonumber
\end{eqnarray}
For $g$ as above, and $\exp z=\left( \begin{array}{cc}
I_n &  z\\
0 & I_m \\
\end{array}\right)$, 

\begin{eqnarray}
g \exp z=\left( \begin{array}{cc}
A & Az+B \\
C & Cz+D \\
\end{array}\right),
\end{eqnarray}
and hence 
\begin{eqnarray}
(g\exp z)_+
=\left( \begin{array}{cc}
I_n & (Az+B)(Cz+D)^{-1} \\
0 & I_m\\
\end{array}\right)
\end{eqnarray}
by \eqref{konkret-HC}.

We also use the Harish-Chandra decomposition to describe the differentials $dg(z)$ for group elements $g$ at points $z$.
We identify all tangent spaces $T_z(\mathscr{D})$ with $\mathfrak{p}^+(\cong M_{nm}(\mathbb{C}))$. Then 
$dg(z): \mathfrak{p}^+ \rightarrow \mathfrak{p}^+$ is
given by the mapping
\begin{eqnarray}
dg(z)=\mbox{Ad}((g\exp z)_0)|_{\mathfrak{p}^+}
\end{eqnarray}
(cf. \cite{satake}).
In the explicit terms given by \eqref{konkret-HC}, this mapping is given by
\begin{eqnarray*}
dg(z)Y=(A-(Az+B)(Cz+D)^{-1}C)YD^{-1},\quad Y \in M_{nm}(\mathbb{C}).
\end{eqnarray*}

\subsection{Strongly orthogonal roots}
We recall that two roots, $\alpha$ and $\beta$, are \emph{strongly orthogonal} if neither $\alpha+\beta$, nor $\alpha-\beta$ is a root.
We define a maximal set of strongly orthogonal noncompact roots, $\Gamma$, inductively by choosing $\gamma_{k+1}$ as the smallest
noncompact root strongly orthogonal to each of the members $\{\gamma_1, \ldots, \gamma_k\}$ already chosen.
When the ordering of the roots in given as in \eqref{ordn}, we get
\begin{eqnarray}
\Gamma=\{\gamma_1, \ldots, \gamma_m\},\, \gamma_j=E_{jj}^*-E_{j+n\,j+n}^*.
\end{eqnarray}
We now let $E_{\gamma_j}$ denote the elementary matrix that spans the root space $\mathfrak{g}^{\gamma_j}$.    
Then the real vector space
\begin{eqnarray}
\mathfrak{a}:=\sum_{j=1}^n \mathbb{R}(E_{\gamma_j}-\theta E_{\gamma_j}) \label{abelsk}
\end{eqnarray}
is a maximal abelian subspace of $\mathfrak{p}$. We set
\begin{eqnarray}
E_j:=E_{\gamma_j}-\theta E_{\gamma_j}. \label{abas}
\end{eqnarray}
\subsection{Shilov boundary}
Let $\mathcal{O}(\mathscr{D})$ denote the set of holomorphic functions on $\mathscr{D}$, and let $\mathcal{O}(\overline{\mathscr{D}})$
denote the subset consisting of those which have continuous extensions to the boundary.
The Shilov boundary of $\mathscr{D}$ is the set
$$\mathcal{S}=\{z \in V | I_m-z^*z=0\}.$$
It has the property that
\begin{eqnarray}
\sup_{z \in \overline{\mathscr{D}}} |f(z)|=\sup_{z \in \mathcal{S}}|f(z)|, f \in \mathcal{O}(\overline{\mathscr{D}}),
\end{eqnarray}
and it is minimal with respect to this property, i.e., no proper subset of $\mathcal{S}$ has the property.
The set $\mathcal{S}$ can also be described as the set of all rank $m$ partial isometries from $\mathbb{C}^m$ to $\mathbb{C}^n$.
The group $K=U(n) \times U(m)$ acts transitively on $\mathcal{S}$ by 
$$(g,h)(z)=gzh^{-1}.$$
To find the isotropy group of the fixed element $z_0:=\left(\begin{array}{c}I_m \\ 0 \end{array}\right)$, 
let $(g,h) \in U(n) \times U(m)$ and write $g$ in the form
$$g=\left(\begin{array}{cc}
A & B\\
C & D\\  
\end{array}\right),$$
where $A$ is of size $m \times m$.
Then
\begin{eqnarray*}
gz_0h^{-1}=\left(\begin{array}{cc}
A & B\\
C & D\\  
\end{array}\right)
\left(\begin{array}{c}h^{-1} \\ 0 
\end{array}\right)
=\left(\begin{array}{c}
Ah^{-1} \\ 
Ch^{-1}
\end{array}\right).
\end{eqnarray*}
So, the equality $gz_0h^{-1}=z_0$ holds if and only if 
$A=h$ and $C=0$. Since $g$ is unitary, the last condition implies that also $B=0$ and hence the isotropy group is
$$K_0:=(U(n) \times U(m))_{z_0}=\left\{(g,h)\in U(n) \times U(m) |
g=\left(\begin{array}{cc}
h & 0\\
0 & D\\  
\end{array}\right)\right\}.$$
Thus we have the description
$$\mathcal{S}=K/K_0=(U(n) \times U(m))/U(n-m)\times U(m)$$
of the Shilov boundary as a homogeneous space.

In the sequel, we will often be concerned with the submanifold $\mathcal{S}_{\Delta}$ of $\mathcal{S}$, where
\begin{eqnarray}
\mathcal{S}_{\Delta}:=\left\{ z_{\underline{\xi}}:=\left(
\begin{array}{ccc}
\xi_1 & \cdots & 0\\
\vdots & \ddots & 0\\
0 & \cdots & \xi_m\\
0 & \cdots & 0\\
\vdots & \vdots & \vdots\\
0 & \cdots & 0\\
\end{array}
\right) \in \mathcal{S}| \xi_1, \ldots, \xi_m \in S^1 \right\}. \label{shilovdiagonal}
\end{eqnarray}
Also, we let $\mbox{diag}(\underline{\xi})$ denote the $m \times m$-matrix
$\left(\begin{array}{ccc}
\xi_1 & \cdots & 0\\
\vdots & \ddots & 0\\
0 & \cdots & \xi_m\\
\end{array}\right)$. 
The identity
\begin{eqnarray}
z_{\underline{\xi}}=
\left(\begin{array}{c}
\mbox{diag}(\underline{\xi})\\
0 \\
\end{array}\right)
=
\left(\begin{array}{cc}
\mbox{diag}(\underline{\xi}) & 0\\
0& I_{n-m}\\
\end{array}\right)
\left(\begin{array}{c}
I_m\\
0 \\
\end{array}\right)
\end{eqnarray}
identifies the matrices in the submanifold $\mathcal{S}_{\Delta}$ with certain cosets in 
$K/K_0$.

\subsection{The real form $\mathscr{X}$}
Consider the mapping $\tau: \mathscr{D} \rightarrow \mathscr{D}$ defined by
\begin{eqnarray}
\tau(z)=\overline{z},
\end{eqnarray}
where the conjugation is entrywise.
It is an antiholomorphic involutive diffeomorphism of $\mathscr{D}$. 
We let $\mathscr{X}$ denote the set of fixed points of $\tau$, i.e., 
\begin{eqnarray}
\mathscr{X}=\{z \in \mathscr{D}|\tau(z)=z\}.
\end{eqnarray}
Moreover, $\tau$ defines an involution, which we also denote by 
$\tau$, of $G$ given by 
\begin{eqnarray}
\tau(g)=\tau g \tau^{-1}.
\end{eqnarray}
We let $H$ denote the set of fixed points, i.e., 
\begin{eqnarray}
H=G^{\tau}=\left\{ g \in G| \tau(g)=g\right\}.
\end{eqnarray}
Clearly, $H=SO(n,m)$, i.e., the elements in $G$ with real entries.
The group $H$ acts transitively on $\mathscr{X}$, and the isotropy group of $0$ in $H$ is $L:=H \cap K$. 
Hence
\begin{eqnarray}
\mathscr{X}\cong H/L.
\end{eqnarray}

\subsection{Minimal representation $\mathscr{H}_1$}
We recall that the Bergman kernel of $\mathscr{D}$ is given by 
\begin{eqnarray}
K(z,w)=\det(I_n-zw^*)^{-(n+m)}.
\end{eqnarray}
It has the transformation property
\begin{eqnarray}
K(gz,gw)=J_g(z)^{-1}K(z,w)\overline{J_g(w)}^{-1},
\end{eqnarray}
where $J_g(z)$ denotes the complex Jacobian of $g$ at $z$.
We let $h(z,w)$ denote the function 
\begin{eqnarray}
h(z,w)=\det(I_n-zw^*).
\end{eqnarray}
Then, for real $\nu$, the kernel 
\begin{eqnarray}
h(\cdot,\cdot)^{-\nu}
\end{eqnarray}
is positive definite if and only if $\nu$ belongs to the Wallach set, $\mathcal{W}$. Here,
\begin{eqnarray}
\mathcal{W}=\{0, 1, \ldots, m-1\} \bigcup (m-1,\infty)
\end{eqnarray}
(cf. \cite{FK}).
The kernel $h(\cdot,\cdot)^{-\nu}$ satisfies the transformation rule
\begin{eqnarray}
h(gz,gw)^{-\nu}=J_{g}(z)^{-\frac{\nu}{n+m}}h(z,w)^{-\nu}\overline{J_g(w)}^{-\frac{\nu}{n+m}}. \label{kernel-transform}
\end{eqnarray}
For $\nu \in \mathcal{W}$, we denote the Hilbert space defined by the kernel $h(\cdot,\cdot)^{-\nu}$ by $\mathscr{H}_{\nu}$.
A projective representation, $\pi_{\nu}$, of $G$ is defined on $\mathscr{H}_{\nu}$ by
\begin{eqnarray}
\pi_{\nu}(g)f(z)=J_{g^{-1}}(z)^{\frac{\nu}{n+m}}f(g^{-1}z).
\end{eqnarray}
We will be concerned with the so called minimal representation, i.e, with the representation $\pi_1$ on the space $\mathscr{H}_1$.

\section{The L-invariants}
For any $\nu \in \mathcal{W}$, let $$\mathscr{H}_{\nu}=\bigoplus_{\underline{k}:=-(k_1 \gamma_1+ \cdots + k_m \gamma_m)} 
\mathcal{P}^{\underline{k}}$$ be
the decomposition into $K$-types.
Here  
$\Gamma=\{\gamma_1, \ldots, \gamma_m\}$ is the maximal strongly orthogonal 
set in $\Delta_{\mathfrak{p}}^+$ with ordering $\gamma_1 < \cdots < \gamma_m$ defined in the previous section, and
\begin{eqnarray}
k_1 \geq \cdots \geq k_m,
k_i \in \mathbb{N},
\end{eqnarray}
and $\mathcal{P}^{\underline{k}}$ is a representation space for the $K$-representation of highest weight that is realised
inside the space of homogeneous polynomials of degree $|\underline{k}|=k_1  + \cdots + k_m $ on $\mathfrak{p}^+$.
When $\nu=1$, the
weights occurring in this sum are all of the form
\begin{eqnarray}
\underline{k}=-k\gamma_1 \label{min-K}
\end{eqnarray}
(cf. \cite{FK}).
Taking $L$-invariants, we have
$$\mathscr{H}_1^L=\bigoplus_{\underline{k}} (\mathcal{P}^{\underline{k}})^L.$$
The data $(K,L, \tau)$ defines a Riemannian symmetric pair, and hence $(V^{\underline{k}})^L$ is at most 
one dimensional by the Cartan-Helgason theorem (cf. \cite{helg2}, Ch. IV, Lemma 3.6.).

We recall the compact Cartan subalgebra $\mathfrak{t} \subset \mathfrak{k}$ in \eqref{cartan}.
We let $\tilde{\mathfrak{t}}$ denote the Cartan subalgebra of $\mathfrak{u}(n) \oplus \mathfrak{u}(m)$ consisting 
of all diagonal imaginary matrices, i.e., matrices of the form \eqref{cartan} but without the requirement that the trace be zero.
Then we have an orthogonal decomposition
\begin{eqnarray}
\tilde{\mathfrak{t}}=\mathfrak{t} \oplus \mathfrak{t}^{\perp}
\end{eqnarray}
given by the Killing form.

Any linear functional $l \in \mathfrak{t}^*$ extends uniquely to a functional on $\tilde{\mathfrak{t}}$ which
annihilates the orthogonal complement $\tilde{\mathfrak{t}}^{\perp}$. We will denote these extensions by the same letter $l$.
Therefore, any dominant integral weight on $\mathfrak{t}$ parametrises an irreducible representation of 
$\mathfrak{u}(n) \oplus \mathfrak{u}(m)$ in which $\tilde{\mathfrak{t}}^{\perp}$ acts trivially.
When $\lambda=\underline{k}=k\gamma_1$ is a $K$-type occurring in $\mathscr{H}_1$, we denote the underlying representation
space for $\mathfrak{u}(n) \oplus \mathfrak{u}(m)$ by $V^{\lambda}$. 
Moreover, the Cartan subalgebra $\tilde{\mathfrak{t}}$ is the sum
$$ \mathfrak{t}=\mathfrak{t}_1 \oplus \mathfrak{t}_2$$ of
the corresponding subalgebras of $\mathfrak {u}(n)$ and  $\mathfrak {u}(m)$ respectively.
The restrictions of $\lambda$ to $\mathfrak{t}_1$ and $\mathfrak{t}_2$ respectively define integral weights, hence
they parametrise irreducible representations of the Lie algebras $\mathfrak{u}(n)$ and $\mathfrak{u}(m)$ respectively.
We denote the corresponding representation spaces by $V_n^{\lambda}$ and $V_m^{\lambda}$.
In what follows, $\lambda$ will always denote the extension to $\mathfrak{u}(n) \oplus \mathfrak{u}(m)$ of a weight  
of the form $\underline{k}$ in \eqref{min-K}. We will use the explicit realisations
\begin{eqnarray}
V_n^{\lambda}=\stackrel{k}{\odot} \mathbb{C}^n,
\end{eqnarray}
where the right hand side denotes the symmetric tensor product defined as a quotient of the $k$-fold tensor product of $\mathbb{C}^n$.
In the following, for a multiindex $\alpha=(\alpha_1,\cdots,\alpha_n)\in \mathbb{N}^n$, we let
\begin{eqnarray}
|\alpha|&:=&\alpha_1+\cdots+\alpha_n,\\
\alpha!&:=&\alpha_1!\cdots\alpha_n!.
\end{eqnarray}
For any choice of orthonormal basis $\{e_1, \ldots, e_n \}$ for $\mathbb{C}^n$, the set
\begin{eqnarray}
\{e^{\alpha}:=e_1^{\alpha_1}\cdots e_n^{\alpha_n}| |\alpha|=k\}
\end{eqnarray}
furnishes a basis for $\stackrel{k}{\odot} \mathbb{C}^n$. We fix an $K$-invariant inner product, $\| \cdot \|_{\mathcal{F}}$, 
\footnote{This is often called the \emph{Fock-Fischer} inner product (cf. \cite{FK}).}
on $\stackrel{k}{\odot} \mathbb{C}^n$
by the normalisation
\begin{eqnarray}
\|e_1^k\|_{\mathcal{F}}^2=k!. 
\end{eqnarray}
Observe that we have suppressed both the indices $k$ and $n$ here. For $n$ fixed, the norm in fact equals the restriction of the
norm defined on all polynomial functions on $\mathbb{C}^n$ (we use the natural identification $e^{\alpha} \leftrightarrow z^{\alpha}$ of
symmetric tensor power with polynomial functions) 
\begin{eqnarray}
\langle p, q \rangle_k:=p(\partial)(q^*)(0), \label{fock}
\end{eqnarray}
where $p(\partial)$ is the differential operator defined by substituting $\frac{\partial}{\partial e_j}$ for $e_j$ in $p$, and 
for $q=\sum_{\alpha}a_{\alpha}z^{\alpha}$, $q^*$ is defined as 
\begin{eqnarray}
(\sum_{\alpha}a_{\alpha}z^{\alpha})^*:=\sum_{\alpha} \overline{a_{\alpha}}z^{\alpha}.
\end{eqnarray}
The suppressing of the index $n$ will not cause any confusion in what follows.
Finally, on the dual space $V_m^{\lambda}$ we have the corresponding basis
\begin{eqnarray}
\{(e^*)^{\alpha}:=(e_1^*)^{\alpha_1}\cdots (e_n^*)^{\alpha_n}| |\alpha|=k\},
\end{eqnarray}
where $\{e_1^*, \ldots, e_n^*\}$ is the dual basis to $\{e_1, \ldots, e_n\}$ with respect to the standard inner product on
$\mathbb{C}^n$. We also let $\| \cdot\|_{\mathcal{F}}$ denote the $K$-invariant norm on $V_m^{\lambda}$ normalised by
\begin{eqnarray}
\|(e_1^*)^k\|_{\mathcal{F}}^2=k!.
\end{eqnarray}

\begin{lemma}
For any choice of orthonormal basis $\{e_1, \ldots, e_m\}$ for $\mathbb{C}^m$ and 
extension $\{e_1, \ldots, e_m, e_{m+1}, \ldots ,e_n\}$ to an orthonormal basis for 
$\mathbb{C}^n$, the vector 
$$\iota_{\lambda}:=\sum_{\stackrel{\alpha \in \mathbb{N}^m}{|\alpha|=k}}
f_{\alpha} \otimes f^*_{\alpha} \in V_n^{\lambda} \otimes (V_m^{\lambda})^*,$$
where $f_{\alpha}=\frac{e^{\alpha}}{(\alpha !)^{1/2}}$ and $f^*_{\alpha}=\frac{(e*)^{\alpha}}{(\alpha !)^{1/2}},$
is $K_0$-invariant.
\end{lemma}

\begin{proof}
We recall the identification of the isotropic subgroup of the fixed element $z_0$ with $U(n-m) \times U(m)$. 
From this it is clear that it suffices to prove that the vector $\iota_{\lambda} \in V_m^{\lambda} \otimes (V_m^{\lambda})^*
\subset V_n^{\lambda} \otimes (V_m^{\lambda})^*$ is invariant under the restriction of the representation of $U(m) \times U(m)$ to
the diagonal subgroup.

The vector space $V_m^{\lambda} \otimes (V_m^{\lambda})^*$ is naturally isomorphic to $\mbox{End}(V_m^{\lambda})$, the isomorphism being
given by $(u \otimes v^*) (y)=v^*(y)u$. Then, if $y \in V_m^{\lambda}$ is the linear combination $y=\sum_{\beta}c_{\beta}f_{\beta}$,
$$\sum_{\alpha}f_{\alpha} \otimes f^*_{\alpha}(y)=\sum_{\alpha,\beta}c_{\beta}\langle f_{\beta},f_{\alpha} \rangle f_{\beta}=y;$$
i.e., $\iota_{\lambda}$ corresponds to the identity operator.
Moreover, for the action of $\mathfrak{u}(m)$ on the tensor product $V_m^{\lambda} \otimes (V_m^{\lambda})^*$, we have
\begin{eqnarray*}
X(u \otimes v^*)(y)&=&(Xu \otimes v^*)(y)+(u \otimes Xv^*)(y)\\
&=&v(y)Xu+(Xv^*)(y)u\\
&=&v(y)Xu+ \langle y ,Xv\rangle u\\
&=&v(y)Xu-\langle Xy ,v\rangle u\\
&=&[X,u \otimes v^*](y),
\end{eqnarray*}
where $X \in \mathfrak{u}(m), u \in V_m^{\lambda}, v^* \in (V_m^{\lambda})^*,$ i.e., the action as derivations of 
the tensor product corresponds to the commutator action on the endomorphisms. In particular, $X \iota_{\lambda}=0$ 
for all $X$ in $\mathfrak{u}(m)$. This proves the lemma.
\end{proof}
Since the vectors in the representation space $V^{\lambda}$ are holomorphic polynomials, they are determined by their restrictions 
to the Shilov boundary $\mathcal{S}$.


In the sequel, we use the Fock inner product to define an antilinear identification of $V_m^{\lambda}$ with $(V_m^{\lambda})^*$  by
\begin{eqnarray*}
v \mapsto v^*, \quad v^*(w)=\langle w,v \rangle_{\mathcal{F}}, \qquad w \in  V_m^{\lambda}.
\end{eqnarray*}
We let $\langle \cdot, \cdot \rangle$ denote the inner product on the tensor product $V_n^{\lambda}\otimes (V_m^{\lambda})^*$
induced by the Fock inner products on the factors. 

\begin{prop}
The operator $T_{\lambda}: V_n^{\lambda} \otimes (V_m^{\lambda})^* \rightarrow V^{\lambda}$ defined by
\begin{eqnarray*}
T_{\lambda}(u \otimes v^*)(z)=\langle (g,h)\iota_{\lambda}, u \otimes v^*\rangle,
\end{eqnarray*}
where $z=(g,h)K_0 \in \mathcal{S}$,
is a $\mathbb{C}$-antilinear isomorphism of $U(n) \times U(m)$-representations. 
\end{prop}

\begin{proof}
We first observe that the left hand side is well defined as a function of $z$ by the invariance of $\iota_{\lambda}$.

The root system $\Delta(\mathfrak{u}(n) \oplus \mathfrak{u}(m), \mathfrak{t})$ is the union of the root systems 
$\Delta(\mathfrak{u}(n), \mathfrak{t}_1)$ and
$\Delta(\mathfrak{u}(m), \mathfrak{t}_2)$. Fix choices of positive roots $\Delta^+(\mathfrak{u}(n), \mathfrak{t}_1)$, and 
$\Delta^+(\mathfrak{u}(m), \mathfrak{t}_2)$ respectively. We define a system of positive roots in 
$\Delta(\mathfrak{u}(n) \oplus \mathfrak{u}(m), \mathfrak{t})$ by
$$\Delta^+(\mathfrak{u}(n) \oplus \mathfrak{u}(m), \mathfrak{t}):=\Delta^+(\mathfrak{u}(n), \mathfrak{t}_1) \cup
\Delta^+(\mathfrak{u}(m), \mathfrak{t}_2).$$
Let $u_{\lambda} \in V_n^{\lambda}$ be a lowest weight-vector, and  $v_{\lambda} \in V_m^{\lambda}$ be
a highest weight-vector. Then $u_{\lambda} \otimes v_{\lambda}^*$ is a 
lowest weight-vector 
in $V_n^{\lambda} \otimes (V_m^{\lambda})^*$. For $H=(H_1,H_2) \in \mathfrak{t}_1 \oplus \mathfrak{t}_2$ we have
\begin{eqnarray*}
&&\frac{d}{dt}(T_{\lambda}(u_{\lambda} \otimes v_{\lambda}^*))(\exp tH \cdot z)_{t=0}\\
&=&\frac{d}{dt}\langle(\exp tH_1g,\exp tH_2h)\iota_{\lambda}, u_{\lambda} \otimes v_{\lambda}^*\rangle_{t=0}\\
&=&\frac{d}{dt}\left(\langle (g,h)\iota_{\lambda}, (\exp -tH_1, \exp -tH_2)(u_{\lambda} \otimes v_{\lambda}^*) \rangle\right)_{t=0}\\
&=& \langle (g,h)\iota_{\lambda}, \lambda(-H_1)u_{\lambda} \otimes v_{\lambda}^* \rangle\\
&&+ \langle (g,h)\iota_{\lambda}, u_{\lambda} \otimes \lambda(-H_2)v_{\lambda}^* \rangle\\
&=&\lambda(H)T_{\lambda}(u_{\lambda} \otimes v_{\lambda}^*)(z).
\end{eqnarray*}
Thus $T_{\lambda}(u_{\lambda} \otimes v^*_{\lambda})$ is a vector of weight $\lambda$. 

Any root vector in $\mathfrak {u}(n) \oplus \mathfrak {u}(m)$ lies in either of the components. Take therefore 
a positive root vector $E+iF \in \mathfrak{u}(n)^{\mathbb{C}}$. Then
\begin{eqnarray*}
&&(E+iF,0)(T_{\lambda}(u_{\lambda} \otimes v_{\lambda}^*))(z)\\
&=&\frac{d}{dt}\left(\langle (\exp t Eg,h)\iota_{\lambda}, u_{\lambda} \otimes v_{\lambda}^* \rangle\right)_{t=0}\\
&&+i\frac{d}{dt}\left(\langle (\exp t Fg,h)\iota_{\lambda}, u_{\lambda} \otimes v_{\lambda}^* \rangle\right)_{t=0}\\
&=&\langle (g,h)\iota_{\lambda}, (-(E-iF)u_{\lambda} \otimes v_{\lambda}^*) \rangle\\
&=&0,
\end{eqnarray*}
since $E-iF$ is a negative root vector.
Similarly one shows that the positive root vectors in  $\mathfrak {u}(m)$ annihilate  $T_{\lambda}(u_{\lambda} \otimes v_{\lambda}^*)$.
The function $T_{\lambda}(u_{\lambda} \otimes v_{\lambda}^*)$ on the Shilov boundary naturally extends to a holomorphic polynomial
on $\mathscr{D}$ which belongs to $\mathscr{H}_1$.
Hence  $T_{\lambda}(u_{\lambda} \otimes v_{\lambda}^*)$ can be written as finite sum of highest weight-vectors from the 
$K$-types of $\mathscr{H}_1$. But it is a vector of weight 
$\lambda$, and so by the multiplicity-freeness of the $K$-type decomposition, 
$T_{\lambda}(u_{\lambda} \otimes v_{\lambda}^*)$ is a highest weight-vector in 
$V^{\lambda}$. 
\end{proof}

\begin{lemma}
The space $(V^{\lambda})^L$ is nonzero if and only if $\lambda=-2k\gamma_1$ for $k \in \mathbb{N}$.
In this case, it is one-dimensional with a basis vector $\psi_k$, where
\begin{eqnarray}
\psi_k(z_{\underline{\xi}}):=\sum_{\stackrel{\beta \in \mathbb{N}^m}{|\beta|=k}} \left( 
\begin{array}{c}
k\\
\beta\\
\end{array}\right)^2 (2\beta!)\xi^{2\beta}, \label{basvektor_k}
\end{eqnarray}
where $z_{\underline{\xi}}$ is the matrix defined in \eqref{shilovdiagonal}.
\end{lemma}

\begin{proof}
We use the isomorphism from the proposition above. Then the first statement is obvious, since for any $\lambda=-j\gamma_1$, the
representation space $V_n^{\lambda}$ is isomorphic to the space of all polynomials of homogeneous degree $j$ on $\mathbb{C}^n$, 
and the corresponding statement holds for $V_m^{\lambda}$. Assume therefore that $\lambda=-2k\gamma_1$.

Clearly, the vector $(e_1^2+ \cdots +e_n^2)^{k} \otimes ((e_1^*)^2+ \cdots +(e_m^*)^2)^{k}$ is an $L$-invariant vector in
$V_n^{\lambda} \otimes (V_m^{\lambda})^*$. We compute its image under $T_{\lambda}$ when restricted to the matrices in 
$\mathcal{S}_{\Delta}$.
\begin{eqnarray*}
T_{\lambda}((e_1^2&+&\cdots +e_n^2)^{k} \otimes ((e_1^*)^2+ \cdots +(e_m^*)^2)^{k})(z_{\underline{\xi}} )\\
&=&\langle (g_{\underline{\xi}}, I_m)\iota_{\lambda}, 
(e_1^2+ \cdots +e_n^2)^{k} \otimes ((e_1^*)^2+ \cdots +(e_m^*)^2)^{k} \rangle\\
&=&\langle \sum_{\alpha}\xi^{\alpha}f_{\alpha} \otimes f^*_{\alpha}, 
(e_1^2+ \cdots +e_n^2)^{k} \otimes ((e_1^*)^2+ \cdots +(e_m^*)^2)^{k} \rangle\\
&=&\sum_{\alpha}\xi^{\alpha}\langle f_{\alpha},(e_1^2+ \cdots +e_n^2)^{k}\rangle \langle f^*_{\alpha}, 
((e_1^*)^2+ \cdots +(e_m^*)^2)^{k} \rangle.
\end{eqnarray*}
Since the symmetric tensor $(e_1^2+ \cdots +e_n^2)^{k}$ has the monomial expansion
$$(e_1^2+ \cdots +e_n^2)^{k}=\sum_{|\beta|=k} \left( 
\begin{array}{c}
k\\
\beta\\
\end{array}\right)e^{2\beta},$$ 
we get the equality
$$T_{\lambda}((e_1^2+ \cdots +e_n^2)^{k} \otimes ((e_1^*)^2+ \cdots +(e_m^*)^2)^{k})(z_{\underline{\xi}})=
\sum_{|\beta|=k} \left( 
\begin{array}{c}
k\\
\beta\\
\end{array}\right)^2 (2\beta!)\xi^{2\beta}.$$
\end{proof}

\begin{thm}
The polynomials $\varphi_{k}$ of degree $2k$, for $k \in \mathbb{N}$, given by
$$\varphi_k(z_{\underline{\xi}})=\frac{1}{4^k k!\left(\frac{m}{2}\right)_k^{1/2}\left(\frac{n}{2}\right)_k^{1/2}}
\sum_{|\beta|=k} \left( 
\begin{array}{c}
k\\
\beta\\
\end{array}\right)^2 (2\beta!)\xi^{2\beta}$$
constitute an orthonormal basis for the subspace, $\mathscr{H}_1^L$,  of $L$-invariants.
\end{thm}

\begin{proof}
The only thing that is left to prove is the normalisation part of the statement, i.e., we need to compute the norms of the polynomials
$\psi_k$.

Using the antilinear isomorphism $T_{\lambda}$, we can introduce an inner product 
\begin{eqnarray*}
\langle \cdot, \cdot  \rangle^{'}_{\lambda}:=\overline{\langle T^{-1}_{\lambda} \cdot, T^{-1}_\lambda \cdot \rangle},
\end{eqnarray*}
where the
the right hand side denotes the conjugate of the inner product on the tensor product induced by the Fock inner products 
on the factors, on $V^{\lambda}$.
By Schur's lemma, the equality
$$\| \cdot \|_{\mathcal{F}}=C_{\lambda}\| \cdot \|^{'}_{\lambda}$$ holds on $V^{\lambda}$ for some complex constant $C_{\lambda}$.
To compute this constant, we compare the norms of the lowest weight-vector $u_{\lambda} \otimes v_{\lambda}^*$ and 
the highest weigh-vector $T_{\lambda}(u_{\lambda} \otimes v_{\lambda}^*)$ in their respective representation spaces.
Let $\{e_1, \ldots, e_m\}$ and $\{e_1, \ldots, e_n\}$ denote the standard orthonormal bases 
for $\mathbb{C}^n$ and $\mathbb{C}^m$
respectively. Then $u_{\lambda} \otimes v^*_{\lambda}=e_1^{2k} \otimes (e_1^{2k})^*$, and
\begin{eqnarray*}
\| e_1^{2k} \otimes (e_1^{2k})^*\|=(2k) !.
\end{eqnarray*}
Moreover, the normalised lowest weight-vector $\frac{e^{2k}_1 \otimes (e^*_1)^{2k}}{(2k) !}$ maps to
\begin{eqnarray*}
T_{\lambda}\left(\frac{e^{2k}_1 \otimes (e^*_1)^{2k}}{(2k) !}\right),
\end{eqnarray*}
where
\begin{eqnarray*}
T_{\lambda}\left(\frac{e^{2k}_1 \otimes (e^*_1)^{2k}}{(2k) !}\right)(z_{\underline{\xi}})&=&\xi_1^{2k}\\
&=&p_{11}(z_{\underline{\xi}}),
\end{eqnarray*}
where $p_{11}$ is the highest weight vector given by $p_{11}(z)=z_{11}^{2k}$. 
Since $\|p_{11}\|_{\mathcal{F}}=\sqrt{(2k)!}$, we see that $C_{\lambda}=\sqrt{(2k)!}$.

The norm of $(e_1^2+\cdots +e_n^2)^k \otimes (((e_1^*)^2+\cdots +(e_m^*)^2)$ is straightforward to compute.
In fact, 
\begin{eqnarray*}
\|(e_1^2+ \cdots +e_n^2)^k\|_{\mathcal{F}}^2 \, \|((e_1^*)^2+\cdots +(e_m^*)^2)^k\|_{\mathcal{F}}^2=
(k!)^2\left(\frac{m}{2}\right)_k\left(\frac{n}{2}\right)_k.
\end{eqnarray*}
Finally, we have the equality
\begin{eqnarray}
\| \cdot \|_{1}^2=\frac{1}{(2k)!}\| \cdot \|_{\mathcal{F}}^2
\end{eqnarray}
(cf. \cite{FK})
relating the $\mathscr{H}_1$-norm to the Fock-Fischer norm on the $K$-type $\underline{2k}=-2k\gamma_1$,
and this ends the proof.
\end{proof}

\section{The action of the Casimir element on the L-invariants}
We consider the representation of the universal enveloping algebra $U(\mathfrak{h}^{\mathbb{C}})$ defined for all $X \in \mathfrak{h}$ by
\begin{eqnarray}
f \mapsto \frac{d}{dt}\pi_1(\exp tX)f|_{t=0},
\end{eqnarray}
for $f$ in the dense subspace, $\mathscr{H}_1^{\infty}$, of analytic vectors, and extended to a homomorphism 
$U(\mathfrak{h}^{\mathbb{C}}) \rightarrow \mbox{End}(\mathscr{H}_1^{\infty})$. We will denote this representation too by $\pi_1$.
We recall that the Casimir element, $\mathcal{C} \in U(\mathfrak{h}^{\mathbb{C}})$ is given by
\begin{eqnarray}
\mathcal{C}=X_1^2+\cdots +X_p^2-Y_1^2-\cdots-Y_q^2,
\end{eqnarray}
 where $\{X_i,i=1, \ldots, \mbox{dim}\mathfrak{q}\}$ and $\{Y_i, i=1, \ldots, \mbox{dim}\mathfrak{l}\}$ are any orthogonal bases 
for $\mathfrak{q}$ and
$\mathfrak{l}$ respectively with respect to the Killing form, $B( \cdot, \cdot)$, on $\mathfrak{h}$ such that
\begin{eqnarray*}
B(X_i,X_i)&=&1,\qquad i=1, \ldots, \mbox{dim}\,\mathfrak{q},\\ 
B(Y_i,Y_i)&=&-1, \qquad i=1, \ldots, \mbox{dim}\,\mathfrak{l}.
\end{eqnarray*}

Consider now the left regular representation, $l$, of $H$ on $C^{\infty}(H/L)$, i.e., $l(h)f(x)=f(h^{-1}x)$. 
We define an operator $R_1: \mathscr{H}_1 \rightarrow C^{\infty}(H/L)$ by
\begin{eqnarray}
R_1f(x):=h(x,x)^{-1/2}f(x). \label{R_1}
\end{eqnarray}
This is the \emph{generalised Segal-Bargmann transform} due to Ólafsson and \O rsted (cf. \cite{ol-or}).
A nice introduction to this transform in a more general context can also be found in Ólafssons overview paper \cite{gestur-ancont}.
The following lemma is an immediate consequence of the transformation rule \eqref{kernel-transform}.
\begin{lemma}\label{R_1-ekvivariant}
The operator $R_1: \mathscr{H}_1 \rightarrow C^{\infty}(H/L)$ is $H$-equivariant. 
\end{lemma}
Moreover, the Casimir element acts on $C^{\infty}(H/L)$ as 
the Laplace-Beltrami operator, $\mathcal{L}$, for the symmetric space $H/L$. 
We recall the "polar coordinate map" (cf. \cite{helg1}, Ch.IX)
\begin{eqnarray}
\phi: L/M \times A^{+} \rightarrow (H/L)',\\
(lM,a) \mapsto laL \nonumber
\end{eqnarray}
Here $(H/L)':=H'/L$, where $H'$ is the set of regular elements in $H$, and $A^{+}=\exp \mathfrak{a}^{+}$, where
\begin{eqnarray}
\mathfrak{a}^{+}=\{t_1E_1+\cdots +t_mE_m| t_i \geq 0, i=1, \ldots, m\}.
\end{eqnarray}
The map $\phi$ is a diffeomorphism onto an open dense set in $H/L$. Hence, any $f \in C^{\infty}(H/L)^L$ is uniquely
determined by its restriction to the submanifold $A^+ \cdot 0=\psi(\{eM\} \times A^+)$. In fact, the restriction mapping
$f \mapsto f|_{A^+ \cdot 0}$ defines an isomorphism between the spaces $C^{\infty}(H/L)^L$ and 
$C^{\infty}(A^+ \cdot 0)^{N_L(\mathfrak{a})/Z_L(\mathfrak{a})}$.
The space $C^{\infty}(H/L)^L$ is invariant under the Laplace-Beltrami operator.
Recall that the radial part of the Laplace-Beltrami operator is a differential operator, $\Delta\mathcal{L}$,
on the 
submanifold $A^+ \cdot 0$ with the property that the diagram 
\begin{eqnarray*}
\begin{CD}
 C^{\infty}(H/L)  @>\mathcal{L}>>  C^{\infty}(H/L)\\
@VVV      @VVV\\
C^{\infty}(A^+ \cdot 0) @>\Delta \mathcal{L}>> C^{\infty}(A^+ \cdot 0)
\end{CD},
\end{eqnarray*}
where the vertical arrows denote the restriction map,
commutes.

Moreover, the functions in $\mathscr{H}_1$ are determined by their restrictions to the real submanifold $H/L$, 
and the $L$-invariant functions are determined by their restrictions to $A^+ \cdot 0$. By Lemma \ref{R_1-ekvivariant} and the
above discussion, we have the following commuting diagram.

\begin{eqnarray*}
\begin{CD}
 \mathscr{H}_1^L  @>\pi_1(C)>> \mathscr{H}_1^L  \\
@VVV      @VVV\\
C^{\infty}(A^+ \cdot 0) @>R_1^{-1}\Delta \mathcal{L}R_1>> C^{\infty}(A^+ \cdot 0)
\end{CD},
\end{eqnarray*}
where, again, the vertical arrows denote the restriction maps.

In what follows, we will compute the action of the operator $R_1^{-1}\Delta \mathcal{L}R_1$ on the subspace $\mathscr{H}_1^L$.

The radial part of the Laplace-Beltrami operator of $H/L$ is given by (cf.\cite{helg2}, Ch. II, Prop. 3.9)

\begin{eqnarray*}
&&4\Delta \mathcal{L}=\sum_{j=1}^m\frac{\partial^2}{\partial t_j^2}+\sum_{m \geq i \geq j \geq 1}\coth (t_i \pm t_j)
(\frac{\partial}{\partial t_i}\pm \frac{\partial}{\partial t_j})\\
&&+(n-m)\sum_{j=1}^m \coth t_j \frac{\partial}{\partial t_j}.
\end{eqnarray*}
The coordinates $t_i$ are related to the Euclidean coordinates $x_i$ by $x_i =\tanh t_i$, i.e., 
\begin{eqnarray}
A^+ \cdot 0=\{(x_1, \ldots, x_m) | 0 \leq x_1 \leq x_2 \leq \cdots \leq x_m < 1\}
\end{eqnarray}
In the coordinates $x_i$, the operator $4R_{1}^{-1}\Delta \mathcal{L}R_{1}:=4\mathcal{L}^1$ has the expression 
\begin{eqnarray*}
4\mathcal{L}^1&=&\sum_{i=1}^m\left( -(1-x_i^2)-x_i^2-2x_i(1-x_i^2)\frac{\partial}{\partial x_i}
+(1-x_i^2)^2\frac{\partial^2}{\partial x_i^2}\right)\\
&&+\sum_{i=1}^m\left(2x_i^2-2x_i(1-x_i^2)\frac{\partial}{\partial x_i} \right)\\
&&+(n-m)\sum_{i=1}^m\left(-1-x_i\frac{\partial}{\partial x_i}+ \frac{1}{x_i}\frac{\partial}{\partial x_i} \right)\\
&&+2\sum_{m \geq i > j \geq 1}\left(-1+\frac{(1-x_i^2)(1-x_j^2)}{x_i^2-x_j^2}\left(x_i\frac{\partial}{\partial x_i}-
x_j\frac{\partial}{\partial x_j}\right)\right).
\end{eqnarray*}
The following lemma is proved by a straightforward calculation. A proof for a similar decomposition can be found in \cite{zhtp}.

\begin{lemma}
The operator $4R_{1}^{-1}\Delta\mathcal{L}R_{1}$ can be written as a sum of three operators, $\mathcal{L}_{-}, \mathcal{L}_{0}$ and 
$\mathcal{L}_{+}$ that lower, keep and, respectively, raise the degrees of the polynomials $\psi_k$. In fact,
\begin{eqnarray*}
\mathcal{L}_{-}&=&\sum_{i=1}^m\left(\frac{\partial^2}{\partial x_i^2}+\frac{n-m}{x_i}\frac{\partial}{\partial x_i}\right)
+2\sum_{m \geq i > j \geq 1}\frac{1}{x_i^2-x_j^2}\left(x_i\frac{\partial}{\partial x_i}-x_j\frac{\partial}{\partial x_j}\right),\\
\mathcal{L}_0&=&-mn+\sum_{i=1}^m \left((-4-(n-m))x_i\frac{\partial}{\partial x_i}-2x_i^2
\frac{\partial^2}{\partial x_i^2}\right)\\
&&-2\sum_{m \geq i > j \geq 1}\frac{x_i^2+x_j^2}{x_i^2-x_j^2}\left(x_j\frac{\partial}{\partial x_j}-x_i\frac{\partial}{\partial x_i}\right),\\
\mathcal{L}_{+}&=&\sum_{i=1}^m\left(2x_i^2+4x_i^3\frac{\partial}{\partial x_i}+x_i^4\frac{\partial^2}{\partial x_i^2}\right)\\
&&+2\sum_{m \geq i > j \geq 1}\frac{x_i^2x_j^2}{x_i^2-x_j^2}\left(x_i\frac{\partial}{\partial x_i}-x_j\frac{\partial}{\partial x_j}\right).
\end{eqnarray*} 
\end{lemma}

\begin{prop}
The operator $\mathcal{L}^1$ acts on the (unnormalised) orthogonal system $\{\psi_k\}$ as the Jacobi operator 
$$\mathcal{L}^1\psi_k=A_k\psi_{k-1}+B_k \psi_k+C_k\psi_{k+1},$$
where 
\begin{eqnarray}
\,\,\,\,\,\,\,\,\,\,\,\,\,\,\,\,\,\,\,
A_k&=&4k^4+(4(m-2)+2(n-m))k^3\\ \label{jacobikonstanter}
&&+((m^2-4m+4)+(n-m)(m-2))k^2,\nonumber \\
B_k&=&-2k^2-\frac{n+m}{2}k-\frac{mn}{4},\nonumber \\ 
C_k&=&\frac{1}{4}. \nonumber
\end{eqnarray}
\end{prop}

\begin{proof}
It follows from the above lemma that the operator is a Jacobi operator. In order to identify the coefficients 
$A_k, B_k,$ and $C_k$, 
we evaluate the polynomials at points $(x_1,0):=(x_1,0, \ldots, 0)$.
Then we have 

\begin{eqnarray*}
\mathcal{L}^+\psi_k ((x_1,0))&=&\left(2x_1^2+4x_1^3\frac{\partial_1}{\partial x_1}
+x_1^4\frac{\partial_1^2}{\partial x_1^2}\right)\psi_k((x_1,0))\\
&=&(2+8k+2k(2k-1))(2k)!x_1^{2k+2}\\
&=&\frac{4k^2+6k+2}{(2k+2)(2k+1)}\,\psi_{k+1}((x_1,0))\\
&=&\psi_{k+1}((x_1,0)),
\end{eqnarray*}
whence $C_k=\frac{1}{4}$.

We now investigate the action of the operators $\frac{x_i \frac{\partial}{\partial x_i}-x_j 
\frac{\partial}{\partial x_j}}{x_i^2-x_j^2}$ that occur in $\mathcal{L}_{-}$ and in $\mathcal{L}_0$.
For $i$ and $j$ fixed, we write the symmetric polynomial $\psi_k$ as a sum (suppressing here the indices $k,i$ and $j$
in order to increase readability)
$$\psi_k=\sum_{c \geq d \geq 0}p_{c,d}(x)(x_i^{2c}x_j^{2d}+x_i^{2d}x_j^{2c}),$$
where the $p_{c,d}$ are symmetric polynomials in the variables other than $x_i$ and $x_j$. The operator then acts on the second factor 
of each term, and 
\begin{eqnarray*}
&&\frac{x_i \frac{\partial}{\partial x_i}-x_j 
\frac{\partial}{\partial x_j}}{x_i^2-x_j^2}(x_i^{2c}x_j^{2d}+x_i^{2d}x_j^{2c})\\
&&=2(c-d)(x_ix_j)^{2d}(x_i^{2(c-d-1)}+ \cdots +x_j^{2(c-d-1)}).
\end{eqnarray*}
Evaluating the right hand side at $(x_1,0)$ (whence $x_i=0$) yields zero unless $d=0$, in which case we get $2cx_j^{2(c-1)}$.
Therefore,
\begin{eqnarray*}
&&\frac{1}{x_i^2-x_j^2}\left(x_i\frac{\partial}{\partial x_i}-x_j\frac{\partial}{\partial x_j}\right)
(\psi_k)((x_1,0))\\
&&=\sum_{c=0}^k p_{c,0}((x_1,0))(2cx_j^{2(c-1)})((x_1,0)).
\end{eqnarray*}
We now consider two separate cases.
\begin{enumerate}
\item If $j=1$, then evaluating the polynomial $p_{c,0}$ at a point $(x_1, 0)$ yields zero unless it is a constant polynomial, 
i.e., unless $c=k$. In this case, $p_{k,0}=(2k)!$.\\
\item If $j \neq 1$, then evaluating $p_{c,0}2cx_j^{2(c-1)}$ at $(x_1, 0)$ gives zero unless $c=1$, in which case we get 
the value 
\begin{eqnarray*}
&&2p_{1,0}(x_1,0)=2\left(\frac{k!}{(k-1)!}\right)^2(2(k-1))!2!x_1^{2k-2}\\
&&=4k^2(2(k-1))!x_1^{2k-2}.
\end{eqnarray*}
\end{enumerate}
Hence, we have 
\begin{eqnarray*}
\sum_{m \geq i > j \geq 1}\frac{1}{x_i^2-x_j^2}\left(x_i\frac{\partial}{\partial x_i}-x_j\frac{\partial}{\partial x_j}\right)
(\psi_k)((x_1,0))\\
=(m-1)2k(2k)!x_1^{2k-2}+
\left( \begin{array}{c}
m-1\\
2 \\
\end{array}\right)
4k^2(2(k-1))!x_1^{2k-2}.
\end{eqnarray*}
From this, we conclude that
\begin{eqnarray*}
&&\mathcal{L}_{-}\psi_k((x_1,0))\\
&=&2\left((m-1)2k(2k)!x_1^{2k-2}+
\left( \begin{array}{c}
m-1\\
2 \\
\end{array}\right)
4k^2(2(k-1))!x_1^{2k-2}\right)x_1^{2k-2}\\
&&+ \left(2k(2k-1)(2k)!+(m-1)4k^2(2(k-1))!\right)x_1^{2k-2}\\
&&+ \left((n-m)2k(2k)! + 4(m-1)(n-m) k^2(2(k-1))!\right) x_1^{2k-2}\\
&&+(4(m^2-4m+4)+4(n-m)(m-2))(2(k-1))!x_1^{2k-2},
\end{eqnarray*}
and hence 
\begin{eqnarray*}
A_k&=&4k^4+(4(m-2)+2(n-m))k^3\\
&&+((m^2-4m+4)+(n-m)(m-2))k^2.
\end{eqnarray*}

Similarly, we see that
\begin{eqnarray*}
\mathcal{L}_0\psi_k((x_1,0))&=&(-mn +(-(n-m)-4)2k-4k(2k-1))(2k)!x_1^{2k}\\
&&-2(m-1)2k(2k)!x_1^{2k}\\
&=&(-8k^2+(-4(m-1)-2(n-m)-4)k-mn)\psi_k((x_1,0)),
\end{eqnarray*}
and hence the value of $B_k$.
\end{proof}

\begin{thm}
The Hilbert space $\mathscr{H}_1^L$ is isometrically isomorphic to the Hilbert space $L^2(\Sigma,\mu)$, where
\begin{eqnarray*}
\Sigma=(0,\infty) \cup \{i(\frac{1}{2}-\frac{n-m}{4}+k)|k \in \mathbb{N}, \frac{1}{2}-\frac{n-m}{4}+k <0\},
\end{eqnarray*}
and $\mu$
is the measure defined by 
\begin{eqnarray}
\int_{\Sigma}fd\mu=
\frac{1}{2\pi}\int_0^{\infty}\left|\frac{\Gamma(a+ix)\Gamma(b+ix)\Gamma(c+ix)}{\Gamma(2ix)}\right|^2
f(x)dx+ \label{hahnort}\\
\frac{\Gamma(a+c)\Gamma(c+b)\Gamma(b-c)\Gamma(a-c)}{\Gamma(-2c)}
\times \sum_{\stackrel{j \in \mathbb{N}}{c+j<0}}\frac{(2c)_j(c+1)_j(c+b)_j(c+a)_j}{(c)_j(c-b+1)_j(c-a+1)_j}(-1)^j \nonumber \\
\times f(-(c+j)^2),\nonumber 
\end{eqnarray} 
where the constants $a,b,$ and $c$ are given by
\begin{eqnarray}
a&=&\frac{m-1}{2}+\frac{n-m}{4},\label{hahnparametrar}\\
b&=&\frac{1}{2}+\frac{n-m}{4},\nonumber\\
c&=&\frac{1}{2}-\frac{n-m}{4}. \nonumber
\end{eqnarray}
Under the isomorphism, the operator $\mathcal{L}^1$ corresponds to the multiplication operator $f \mapsto -(a^2+x^2)f$.
\end{thm}

\begin{proof}
We recall the continuous dual Hahn polynomials, $S_k(x^2;a,b,c)$, (cf. \cite{koekoek}) defined by
\begin{eqnarray}
\frac{S_k(x^2;a,b,c)}{(a+b)_k(a+c)_k}=_3F_2\left( 
\begin{array}{c}
-k,a+ix,a-ix\\
a+b,a+c
\end{array}
|1\right).
\end{eqnarray}
Here, $( \cdot )_k$ denotes the \emph{Pochhammer symbol} defined as
\begin{eqnarray*}
(t)_0&=&1,\\
(t)_k&=&t(t+1)\cdots(t+k-1), \,k \in \mathbb{N}^+.
\end{eqnarray*}
Suppressing the parameters and denoting the left hand side above by $\tilde{S}_k(x^2)$, these polynomials 
satisfy the recurrence relation
\begin{eqnarray}
-(a^2+x^2)\tilde{S}_k(x^2)=A'_k\tilde{S}_{k-1}(x^2)+B'_k\tilde{S}_k(x^2)+C'_k\tilde{S}_{k+1}(x^2), \label{Jacobi-rec}
\end{eqnarray}
where the recursion constants $A'_k, B'_k,$ and $C'_k$ are given by
\begin{eqnarray}
A'_k&=&k(k+b+c-1),\\
C'_k&=&(k+a+b)(k+a+c),\\
B'_k&=&-(A'_k+C'_k).
\end{eqnarray}
Under a renormalisation of the form 
\begin{eqnarray*}
S_k(x^2,a,b,c) \mapsto \alpha_kS_k(x^2,a,b,c):=S_k(x^2,a,b,c)^{\alpha},
\end{eqnarray*} 
where ${\alpha_k}$ is some sequence of complex numbers, the corresponding polynomials $\tilde{S^{\alpha}}_k$ will also satisfy
a recurrence relation of the type in \eqref{Jacobi-rec}, with constants, $A^{\alpha}_k, B^{\alpha}_k, C^{\alpha}_k$, given by
\begin{eqnarray}
A^{\alpha}_k&=&\frac{\alpha_{k}}{\alpha_{k-1}}A'_k,\\
B^{\alpha}_k&=&B'_k,\\
C^{\alpha}_k&=&\frac{\alpha_{k}}{\alpha_{k+1}}C'_k.\label{alfa-rek}
\end{eqnarray}
From this we can see that the product $A'_{k+1}C'_k=A^{\alpha}_{k+1}C^{\alpha}_k$ is invariant. 

Consider now the continuous dual Hahn polynomials with $S_k(x^2;a,b,c)$, with the parameters $a,b,c$ from \eqref{hahnparametrar}. 
These polynomials satisfy the orthogonality relation (cf. \cite{koekoek})
\begin{eqnarray}
&&\frac{1}{2\pi}\int_0^{\infty}\left|\frac{\Gamma(a+ix)\Gamma(b+ix)\Gamma(c+ix)}{\Gamma(2ix)}\right|^2 
S_k(x^2;a,b,c)S_l(x^2;a,b,c)dx \nonumber \\ 
&&+\frac{\Gamma(a+c)\Gamma(c+b)\Gamma(b-c)\Gamma(a-c)}{\Gamma(-2c)} \nonumber \\ 
&&\times \sum_{\stackrel{j \in \mathbb{N}}{c+j<0}}\frac{(2c)_j(c+1)_j(c+b)_j(c+a)_j}{(c)_j(c-b+1)_j(c-a+1)_j}(-1)^j \nonumber  \\ 
&&\times S_k(-(c+j)^2;a,b,c)S_l(-(c+j)^2;a,b,c)\nonumber \\
&&\quad =\Gamma(k+a+b)\Gamma(k+a+c)\Gamma(k+b+c)k!\delta_{kl}. \label{hahnort}
\end{eqnarray}
By a straightforward computation one sees that the corresponding constants 
$A'_k, B'_k,$ and $C'_k$ are related to the Jacobi constants $A_k, B_k$, and $C_k$ in \eqref{jacobikonstanter} 
by
\begin{eqnarray*}
A_{k+1}C_k&=&A'_{k+1}C'_k,\\
B_k&=&B'_k.
\end{eqnarray*}
We can thus use \eqref{alfa-rek} to define a sequence ${\alpha_k}$ recursively in such a way that the resulting polynomials
$\tilde{S}^{\alpha}_k$ satisfy the recurrence relation
\begin{eqnarray}
-(a^2+x^2)\tilde{S}^{\alpha}_k(x^2)=A_k \tilde{S}^{\alpha}_{k-1}(x^2)+B_k\tilde{S}^{\alpha}_k(x^2)+C_k\tilde{S}^{\alpha}_{k+1}(x^2)
\end{eqnarray}
with the same Jacobi constants as the operator $4\mathcal{L}^1$.
More precisely, we set 
\begin{eqnarray}
\alpha_0&:=&\left(\Gamma\left(\frac{m}{2}\right)\Gamma\left(\frac{n}{2}\right)\right)^{-1/2},\\
\alpha_{k+1}&:=&\frac{1}{4}\left(k+\frac{m}{2}\right)\left(k+\frac{n}{2}\right)^{-1}\alpha_k.
\end{eqnarray}
Then $\alpha_k=\left(\Gamma(\frac{m}{2})\Gamma(\frac{n}{2})\right)^{-1/2} 4^k\left(\frac{m}{2}\right)_k\left(\frac{n}{2}\right)_k$,
and hence, by \eqref{hahnort}, we have
\begin{eqnarray}
\|\tilde{S}^{\alpha}_k\|_{L^2}^2&=&4^{2k}(k!)^2\left(\frac{m}{2}\right)_k\left(\frac{n}{2}\right)_k\\
&=&\|\psi_k\|_1^2.
\end{eqnarray}
Therefore, the operator $T_0: \mathscr{H}_1^L \rightarrow L^2(\Sigma, d\mu)$ defined by
\begin{eqnarray}
T_0\psi_k=\tilde{S}^{\alpha}_k
\end{eqnarray}
is a unitary operator which diagonalises the restriction of the operator $\mathcal{L}^1$ to $\mathscr{H}_1^L$.
\end{proof}

\begin{thm} \label{irruppdelning}
For each $x \in \Sigma$, there exists a Hilbert space $\mathscr{H}_x$ and an irreducible unitary spherical representation, $\pi_x$, of
$H$ on $\mathscr{H}_x$ such that
\begin{enumerate}
\item If $v_x \in \mathscr{H}_x$ is the canonical spherical vector, then there is an isometric embedding of Hilbert spaces
$L^2(\Sigma, \mu) \subset \int_{\Sigma}\mathscr{H}_xd\mu(x)$ given by
\begin{eqnarray*}
f \mapsto s_f,
\end{eqnarray*}
where $s_f(x):=f(x)v_x$.
\item The operator $T_0$ extends uniquely to an $H$-intertwining unitary operator 
\begin{eqnarray}
T: (\pi_1,\mathscr{H}_1) \rightarrow \left(\int_{\Sigma}\pi_xd\mu(x),\int_{\Sigma}\mathscr{H}_xd\mu(x)\right). \label{spektrdekomp}
\end{eqnarray}
\end{enumerate}
\end{thm}

\begin{proof}
The Banach algebra $L^1(H)$ equipped with convolution as multiplication carries the structure of a Banach 
$^*$-algebra when the involution  $^*$ is defined
as $f^*(h)=\overline{f(h^{-1})}$. The representation $\pi_1$ of $H$ induces a representation of $L^1(H)$ by
\begin{eqnarray}
\pi(f)=\int_H f(h)\pi_1(h)dh.
\end{eqnarray}
If $L^1(H)^{\#}$ denotes the subalgebra of left and right $L$-invariant $L^1$-functions, the closed $C^*$-algebra generated by
$\pi_1(L^1(H)^{\#})$ and the identity operator is a commutative $C^*$-algebra.
Moreover, the Casimir operator $\pi_1(\mathcal{C})$ commutes with all the operators $\pi_1(f)$ for $f \in L^1(H)^{\#}$.
Hence, (by \cite{akh},Vol. I, Thm 1, p. 77), the diagonalisation of the Casimir operator yields a simultaneous diagonalisation of the whole
commutative algebra $\pi_1(L^1(H)^{\#})$.

For $f \in L^1(H)^{\#}$, we let the function $\tilde{f}: \Sigma \rightarrow \mathbb{C}$ be the multiplier corresponding
to the operator $T\pi_1(f)T^{-1} : L^2(\Sigma, \mu) \rightarrow L^2(\Sigma, \mu)$.
For each $x \in \Sigma$, we let $\lambda_x$ denote the multiplicative functional
\begin{eqnarray}
\lambda_x(f):=\tilde{f}(x),
\end{eqnarray}
which clearly is bounded almost everywhere with respect to $\mu$.
The equality
\begin{eqnarray*}
\langle \pi_1(f)\varphi_0, \varphi_0 \rangle_1=\int_{\Sigma}\lambda_x(f)d\mu(x)
\end{eqnarray*}
holds for $f \in L^1(H)^{\#}$, i.e., the positive functional 
\begin{eqnarray}
\Phi_0(f):=\langle \pi_1(f)\varphi_0, \varphi_0\rangle_1, f \in L^1(H)^{\#}
\end{eqnarray} 
is expressed as an integral of characters. 

By \cite{papper1} (Thm. 10) there exists a direct integral decomposition into unitary spherical irreducible
representations of the form \eqref{spektrdekomp}, and it expresses 
the functional $\Phi_0$ as an integral of characters against the corresponding measure. This measure is supported on the characters
given by positive definite spherical functions. 
By \cite{rudfa} (Thm. 11.32), such an integral expression for bounded positive functionals 
is unique, and hence every character $\lambda_x$ can be
expressed by a positive definite spherical function $\phi_x$ as
\begin{eqnarray*}
\lambda_x(f)=\int_Hf(h)\phi_x(h)dh.
\end{eqnarray*}
The rest now follows from the proof of Thm. 10 in \cite{papper1}.

\end{proof}

\section{A subrepresentation of $\pi_1|_H$}
Recall that the boundary $\partial \mathscr{D}$ is the disjoint union of $m$ $G$-orbits. More specifically, for $j=1, \ldots, m$, let $e_j$ 
denote the
$n \times m$ matrix with $1$ at position $(j,j)$ and all other entries zero. Then 
$$\partial \mathscr{D}=\bigcup_{r=1}^mG(e_1+ \cdots +e_r)$$ and the inclusion
$$\overline{G(e_1+ \cdots +e_{r+1}) } \subseteq G(e_1+ \cdots +e_r)$$ holds for $r=1, \ldots, m-1$. 
The Shilov boundary is the $G$-orbit of the rank $m$ partial isometry $e_1+ \cdots +e_m$. It is also the $K$-orbit of this element.
We consider now the "real part", $Y$, of the Shilov boundary, i.e., 
\begin{eqnarray}
Y:=\mathcal{S} \cap M_{nm}(\mathbb{R}).
\end{eqnarray}
Then $Y$ is the homogeneous space $H/P_0$, where $P_0$ is the maximal parabolic subgroup defined by the 
one dimensional subalgebra 
\begin{eqnarray*}
\mathfrak{a}_0=\mathbb{R}(E_1+\cdots+E_m) 
\end{eqnarray*}
of $\mathfrak{a}$ (cf. \eqref{abas}).
We let $P_0=M_0A_0N_0$ 
be the Langlands decomposition. Then $Y$ can also be described as a homogeneous space
$Y=L/L\cap M_0$.
Consider the one dimensional representation with character
\begin{eqnarray}
l \mapsto |\det \mbox{Ad}^{-1}_{\mathfrak{l}/\mathfrak{l}\cap \mathfrak{m}_0}(l)| \label{Ad-kvot}
\end{eqnarray}
of $L \cap M_0$. The induced representation $\mbox{Ind}_{L \cap M_0}^L(|\det \mbox{Ad}^{-1}_{\mathfrak{l}/\mathfrak{l}\cap \mathfrak{m}_0}|)$ 
is realised on the space of sections of the density bundle of $Y=L/L \cap M_0$. 
The representation \eqref{Ad-kvot} is in fact trivial, and this allows us
to define an $L$-invariant section, $\omega$, by
\begin{eqnarray}
\omega(l (L\cap M_0)):=l_{e (L\cap M_0)}\omega_0, \label{inv-densitet}
\end{eqnarray}
where $\omega_0 \neq 0 \in \mathcal{D}(T_{e (L \cap M_0)})$ is arbitrary, where $\mathcal{D}(T_{e (L \cap M_0)})$ denotes
the vector space of densities on $T_{e (L \cap M_0)}$.   The section $\omega$ then corresponds to 
a constant function $F_{\omega}: L \rightarrow \mathbb{C}$. In the usual way, we will sometimes identify $\omega$ 
with the measure it defines by integration against continuous functions. We then use measure theoretic notation and
write $\int_Y \varphi \,d\omega$ for $\int_Y \varphi \,\omega$. Moreover, we choose $\omega_0$ in \eqref{inv-densitet}
so that this measure is normalised.

Using the identification $\mathfrak{l}/\mathfrak{l}\cap \mathfrak{m}_0 \simeq \mathfrak{h}/\mathfrak{p}_0$, the 
representation \eqref{Ad-kvot} extends to the representation $\delta_0$ of $\mathfrak{p}_0$ given by
\begin{eqnarray}
\delta_0(m_0a_0n_0)=|\det(\mbox{Ad}_{\mathfrak{h}/\mathfrak{p}_0}(m_0a_0n_0)^{-1}|.
\end{eqnarray}
Clearly, $\delta_0(m_0a_0n_0)=e^{2\rho_0(\log a_0)}$, where $\rho_0$ denotes the half sum of the restricted roots.
The action of $H$ as pullbacks (actually, the inverse mapping composed with pullback) 
on densities is equivalent to the left action defined by the representation 
$\mbox{Ind}_{P_0}^H(\delta_0)$. For the extension of the function $F_{\omega}$ to a $P_0$ equivariant
function $H \rightarrow \mathbb{C}$ (which we still denote by $F_{\omega}$), we then have 
\begin{eqnarray}
F_{\omega}(k_0m_0a_0n_0)=e^{-2\rho_0(\log a_0)}F_{\omega}(k_0)=e^{-2\rho_0(\log a_0)}F_{\omega}(e).
\end{eqnarray}
From this, it follows that 
\begin{eqnarray}
h^*\omega(l(L \cap M_0)=e^{-2\rho_0(\log A_0(hl)}\omega(l (L \cap M_0)).
\end{eqnarray}
The action of $H$ on $Y$ can either be described on the coset space $H/P_0$ in terms of the Langlands decomposition for $P_0$,
or in terms of the geometric action on the boundary of $\mathscr{D}$ defined by the Harish-Chandra decomposition.
The next proposition expresses the transformation of $\omega$ under $H$ in terms of the latter description.

\begin{lemma}
The density $\omega$ transforms under the action of $H$ as
\begin{eqnarray}
h^*\omega(v)=J_{h}(v)^{\left(\frac{n-1}{n+m}\right)}\omega(v). \label{denstransform}
\end{eqnarray}
\end{lemma}

\begin{proof}
The idea of the proof is to use the (non-unique) factorisation $H=LM_0A_0N_0$ of $H$. We prove that the group $N_0$ fixes the reference point
$e_1+\cdots +e_m$ and acts with Jacobian equal to one on the tangent space at $e_1+ \cdots +e_m$, and the group elements in $M_0$
have Jacobian equal to one at $e_1+ \cdots +e_m$.
By the chain rule for differentiation, it then suffices to prove the statement for all group elements in $A_0$.
 
In the Langlands decomposition $\mathfrak{p}_{\mbox{min}}=\mathfrak{m} \oplus \mathfrak{a} \oplus \mathfrak{n}$ for the minimal 
parabolic subgroup,
the subalgebra $\mathfrak{n}$ is generated by the restricted root spaces
\begin{eqnarray*}
\bigoplus_{m \geq j > i \geq 1}\mathfrak{h}_{E_j^*+E_i^*}&=&
\left\{X_q=\left(\begin{array}{ccc}
-q & 0 & q\\
0 & 0 & 0\\
-q & 0 & q \\
\end{array} \right)| q^t=-q \right\},\\
\bigoplus_{m \geq j > i \geq 1}\mathfrak{h}_{E_j^*-E_i^*}&=&
\left\{X_u=\left(\begin{array}{ccc}
u^t-u & 0 & u+u^t\\
0 & 0 & 0\\
u+u^t & 0 & u^t-u \\
\end{array} \right)| \mbox{$u$ is upper triang.} \right\},\\
\bigoplus_{j=1}^m \mathfrak{h}_{E_j^*}&=&
\left\{X_z=\left(\begin{array}{ccc}
0 & z^t & 0\\
-z & 0 & z\\
0 & z^t & 0\\
\end{array} \right)\right\},
\end{eqnarray*} 
where the matrices are written in blocks in such a way that the block-rows are of height $m, n-m,$ and $m$ respectively, and
the block-columns are of width $m, n-m,$ and $m$ respectively.

In the Langlands decomposition $\mathfrak{m}_0 \oplus \mathfrak{a}_0 \oplus \mathfrak{n}_0$, 
the centraliser, $\mathfrak{m}_0$ of $\mathfrak{a}_0$ is
the direct sum
$$\mathfrak{m}_0=\mathfrak{m} \oplus \bigoplus _{m \geq j > i \geq 1} \mathfrak{h}_{E_j^*-E_i^*},$$
and 
\begin{eqnarray}
\mathfrak{n}_0=\bigoplus_{m \geq j > i \geq 1} \mathfrak{h}_{E_j^*+E_i^*} \oplus 
\bigoplus_{j=1}^m \mathfrak{h}_{E_j^*}. \label{N_0}
\end{eqnarray}
The matrices $X_q$ and $X_z$ commute, so in order to prove that the elements in $N_0$ have Jacobian equal to one 
at $e_1+ \cdots +e_m$, it 
suffices to consider elements of the form
\begin{eqnarray*}
\exp X_q&=&\left(\begin{array}{ccc}
1-q & 0 & q\\
0 & 1 & 0\\
-q & 0 & 1+q \\
\end{array} \right),\\
\exp X_z&=&\left(\begin{array}{ccc}
1-\frac{z^tz}{2} & z^t & \frac{z^tz}{2}\\
-z & 1 & z\\
-\frac{z^tz}{2} & z^t & 1+\frac{z^tz}{2} \\
\end{array} \right)
\end{eqnarray*}
separately.

We have 
\begin{eqnarray*}
\exp X_q \exp (e_1+ \cdots +e_m)&=&
\left(\begin{array}{ccc}
1-q & 0 & q\\
0 & 1 & 0\\
-q & 0 & 1+q \\
\end{array} \right)
\left(\begin{array}{ccc}
1 & 0 & 1\\
0 & 1 & 0\\
0 & 0 & 1 \\
\end{array} \right)\\
&=&
\left(\begin{array}{ccc}
1-q & 0 & 1\\
0 & 1 & 0\\
-q & 0 & 1\\
\end{array} \right).
\end{eqnarray*}
If we write this matrix in the block form $\left( \begin{array}{cc}
A & B\\
C & D \\ 
\end{array} \right)$, then the $K^{\mathbb{C}}$-component in the Harish-Chandra decomposition is given by
\begin{eqnarray*}
\left( \begin{array}{cc}
A-BD^{-1}C & 0\\
0 & D \\
\end{array} \right)=
I_{n+m},
\end{eqnarray*}
and hence 
\begin{eqnarray}
J_{\exp X_q}(e_1 + \cdots + e_m)=1. \label{N-jac 1}
\end{eqnarray}

Next, we consider the action of $\exp X_z$. We have
\begin{eqnarray*}
\exp X_z \exp(e_1+ \cdots +e_m)=\left( \begin{array}{ccc}
1-\frac{z^tz}{2} & z^t & 1\\
-z & 1 & 0\\
-\frac{z^tz}{2} & z^t & 1 \\
\end{array} \right).
\end{eqnarray*}
Here, the $K^{\mathbb{C}}$-component is given by
\begin{eqnarray*}
K^{\mathbb{C}}(\exp X_z \exp(e_1+\cdots +e_m))=\left( \begin{array}{ccc}
1 & 0 & 0\\
-z & 1 & 0\\
0 & 0 & 1 \\
\end{array}\right).
\end{eqnarray*}
The complex differential of $\exp X_z$ at $e_1+\cdots +e_m$ is then the linear mapping
\begin{eqnarray}
d\exp X_z(e_1+ \cdots +e_m)Y=\left( \begin{array}{cc}
1 & 0\\
z & 1 \\
\end{array}\right)Y,
\end{eqnarray}
where we have identified the tangent spaces with $\mathfrak{p}^+=M_{nm}(\mathbb{C})$.
Clearly, the determinant of this mapping is 
\begin{eqnarray}
\det \left( \begin{array}{cc}
1 & 0\\
z & 1 \\
\end{array}\right)^m=1. \label{N-jac 2}
\end{eqnarray}

Consider now the subgroup $M_0$. Its Lie algebra $\mathfrak{m}_0$ is reductive with Cartan involution given by the restriction of 
$\theta$ and
the corresponding decomposition is $$ \mathfrak{m}_0= \mathfrak{m}_0 \cap \mathfrak{l} \oplus \mathfrak{m}_0 \cap \mathfrak{q}.$$
The abelian subalgebra $\mathfrak{a}$ is included in $\mathfrak{m}_0 \cap \mathfrak{q}$, and therefore (cf. \cite{knapp1}, Prop. 7.29)
\begin{eqnarray}
\mathfrak{m}_0 \cap \mathfrak{q}=\bigcup_{l \in M_0 \cap L} \mbox{Ad}(l)\mathfrak{a}. \label{KAK}
\end{eqnarray}
We now investigate the Jacobians of arbitrary group elements in $A$. For $H=t_1E_1+\cdots+t_mE_m$, 
$$\exp H= 
\left( \begin{array} {ccc}
\Delta (\underline {\cosh t}) & 0 & \Delta(\underline{\sinh t})\\
0 & 1 & 0\\
\Delta(\underline{\sinh t}) & 0 & \Delta(\underline{\cosh t})\\  
\end{array}\right),$$
where $\Delta(\underline{\cosh t})$ denotes the $m \times m$ diagonal matrix with entries $\cosh t_1, \ldots, \cosh t_m$, and the 
other blocks are analogously defined. 
Then
\begin{eqnarray*}
&&\exp(t_1E_1+ \cdots +t_mE_m)\exp(e_1+\cdots+e_m)\\
&&=\left( 
\begin{array}{ccc}
\Delta(\underline{\cosh t}) & 0 & \Delta(\underline{\cosh t+ \sinh t})\\
0 & 1 & 0\\
\Delta(\underline{\sinh t}) & 0 & \Delta(\underline{\cosh t+ \sinh t})\\ 
\end{array}
\right)
\end{eqnarray*}
The $K^{\mathbb{C}}$-component is
\begin{eqnarray*}
&&K^{\mathbb{C}}(\exp(t_1E_1+ \cdots +t_mE_m)\exp(e_1+\cdots+e_m))\\
&&=\left( 
\begin{array}{ccc}
\Delta(\underline{e^{-t}}) & 0 & 1\\
0 & 1 & 0\\
0 & 0 & \Delta(\underline{e^t})\\
\end{array}
\right),
\end{eqnarray*}
so the differential $d(\exp(t_1E_1+ \cdots +t_mE_m))(e_1+ \cdots +e_m)$ is the mapping
\begin{eqnarray*}
\left( \begin{array}{c}
Y_1\\
Y_2\\
\end{array}
\right)
\mapsto 
\left( \begin{array}{c}
\Delta(\underline{e^{-2t}})Y_1\\
Y_2\Delta(\underline{e^{-t}})\\
\end{array}
\right),
\end{eqnarray*}
where $Y_1$ is the upper $m \times m$ block of the $n \times m$ matrix in the tangent space.
Counting the multiplicities of the eigenvalues $e^{-t_j}$, we see that
\begin{eqnarray}
J_{\exp(t_1E_1+\cdots+t_mE_m)}(e_1+ \cdots +e_m)=e^{-(n+m)\sum_{j=1}^m t_j}. \label{A-jac}
\end{eqnarray}
If we write $\mathfrak{a}$ as the orthogonal sum $\mathfrak{a}=\mathfrak{a}_0 \oplus (\mathfrak{a}_0 )^{\perp}$ 
(with respect to the Killing form),
then $(\mathfrak{a}_0)^{\perp}$ consists of those $t_1E_1 + \cdots +t_mE_m$ in $\mathfrak{a}$ for which $\sum_{j=1}^mt_j=0$.
From the identities \eqref{N-jac 1}, \eqref{N-jac 2}, \eqref{KAK}, and \eqref{A-jac} we can thus conclude that
\begin{eqnarray}
J_h(e_1+\cdots +e_m)=J_{A_0(h)}(e_1+\cdots +e_m).
\end{eqnarray}

On the other hand, by \eqref{N_0}, 
\begin{eqnarray}
&&2\rho_0(t(E_1+ \cdots +E_m)) )) \label{rho} \\ 
&&=2\frac{m(m-1)}{2}t+m(n-m)t=m(n-1)t, \nonumber
\end{eqnarray}
so
\begin{eqnarray}
&&e^{-2\rho_0(t(E_1+ \cdots +E_m))} \label{rho-jacobian} \\
&&=(J_{\exp(t(E_1+ \cdots + E_m))}(e_1+ \cdots +e_m))^{\frac{n-1}{n+m}}. \nonumber
\end{eqnarray} 

\end{proof}

In what follows, we will define a Hilbert space of functions on the manifold $Y$.
Hilbert spaces of a similar kind were also considered by Neretin and Olshanski in \cite{ner-ol}. 
One difference is that their spaces were not defined
using a limit procedure (see the next definition below).

We begin by introducing some notation.
For a continuous function, $f$, on $Y$ and $r \in (0,1)$, we define the function $F_r:Y \rightarrow \mathbb{C}$ by
\begin{eqnarray}
F_r(u):=\int_Y f(v)\det(I_n-ruv^t)^{-1}d\omega(v).
\end{eqnarray}
We construct the Hilbert space by requiring that the following space of functions be dense.
\begin{defin}  \label{C_0}
Let $\mathscr{C}_0$ denote the set of all continuous functions\\ 
$f: Y \rightarrow \mathbb{C}$ such that the limit function
\begin{eqnarray*}
 F(u):=\lim_{r \rightarrow 1}F_r(u)
\end{eqnarray*}
exists in the supremum norm.
\end{defin}

On $\mathscr{C}_0$ we define a sesquilinear form $\langle \quad,\quad\rangle_{\mathscr{C}_0}$ by
\begin{eqnarray}
\langle f, g \rangle_{\mathscr{C}_0}=\int_Y f(u)\overline{G(u)}d\omega(u).
\end{eqnarray}
By the Dominated Convergence Theorem, we have
\begin{eqnarray}
&&\int_Y f(u)\overline{G(u)}d\omega(u)  \label{DCT} \\
&&=\lim_{r \rightarrow 1}\int_Y f(u) \int_Y \overline{g(v)}\det(I_n-ruv^t)^{-1}d\omega(v)d\omega(u), \nonumber
\end{eqnarray}
and hence the form $\langle \quad,\quad \rangle_{\mathscr{C}_0}$ is positive semidefinite. 
Let $\mathcal{N}$ denote the space of functions of norm zero, i.e., 
\begin{eqnarray}
\mathcal{N}=\{f \in \mathscr{C}_0 | \langle f, f \rangle_{\mathscr{C}_0}=0\}.
\end{eqnarray}
Then the quotient space $\mathscr{C}_0/\mathcal{N}$ together with the induced sesquilinear form, 
$\widetilde{\langle \quad,\quad \rangle_{\mathscr{C}_0}}$, is 
a pre-Hilbert space.
We define $\mathscr{C}$ to be the
Hilbert space completion of $\mathscr{C}_0$ with respect to $\widetilde{\langle \quad,\quad \rangle_{\mathscr{C}_0}}$. 
We denote the inner product on $\mathscr{C}$ by
$\langle \quad, \quad\rangle_{\mathscr{C}}$.

\begin{prop}
The action $\tau$ of $H$ on $\mathscr{C}_0$ given by
\begin{eqnarray}
\tau(h)f(\eta):=J_{h^{-1}}(\eta)^{\beta}f(h^{-1}\eta), \label{tauverkan} 
\end{eqnarray}
where $\beta=\frac{n-2}{n+m}$,
descends to a unitary representation of $H$ on $\mathscr{C}$.
\end{prop}

\begin{proof}
It suffices to prove that the dense subspace $\mathscr{C}_0/\mathcal{N}$ of $\mathscr{C}$ is $H$-invariant and that the 
action is unitary on $\mathscr{C}_0/\mathcal{N}$.
For this, it clearly suffices to prove that the space $\mathscr{C}_0$ is $H$-invariant, and that $H$ preserves the
sesquilinear form $\langle \quad, \quad \rangle_{\mathscr{C}_0}$, since then the subspace $\mathcal{N}$ is also
$H$-invariant.

Consider first the mapping $f \mapsto F$  in Definition \ref{C_0}. We write $K_1$ for the reproducing kernel. For $h \in H$, we 
then have
\begin{eqnarray*}
\int_Y \tau(h)f(v) K_1(ru,v)d\omega(v)&=&\int_Y J_h(h^{-1}v)^{-\beta} f(h^{-1})K_1(ru_1,v)d\omega(v)\\
&=&\int_Y J_h(v')^{-\beta+\frac{n-1}{n+m}} f(v')K_1(ru_1,hv')d\omega(v'),\\
\end{eqnarray*}
by the transformation property for the measure $\omega$. By the transformation rule for the reproducing kernel $K_1$, we have
\begin{eqnarray*}
&&\int_Y J_h(v')^{-\beta+\frac{n-1}{n+m}} f(v')K_1(ru_1,hv')d\omega(v')\\
&&=\int_Y J_h(h^{-1}ru)^{-\frac{1}{n+m}}f(v')K_1(h^{-1}ru,v')d\omega(v').
\end{eqnarray*}
Therefore,
\begin{eqnarray*}
\lim_{r \rightarrow 1}\int_Y \tau(h)f(v)K_1(ru,v)d\omega(v)= J_h(h^{-1}u)^{-\frac{1}{n+m}}F(h^{-1}u),
\end{eqnarray*}
where the convergence is uniform in $u$, so $\mathscr{C}_0$ is $H$-invariant. 

Next, take $f, g \in \mathscr{C}_0$.
Then, $\langle \tau(h)f,\tau(h)g \rangle_{\mathscr{C}_0}$ is given by
\begin{eqnarray*}
\langle \tau(h)f,\tau(h)g \rangle_{\mathscr{C}_0}&=&\int_Y J_h(h^{-1}u) f(h^{-1}u)J_h(h^{-1}u)^{-\frac{1}{n+m}}\overline{G(h^{-1}u)}d\omega(u)\\
&=&\int_Y f(v')\overline{G(v')}d\omega(v')\\
&=&\langle f, g \rangle_{\mathscr{C}_0}, 
\end{eqnarray*}
where the second equality follows from the transformation property of $\omega$.
\end{proof}

The next proposition gives a sufficient condition for the Hilbert space $\mathscr{C}$ to be nonzero.

\begin{prop}
The (equivalence class modulo $\mathcal{N}$ of the) constant function $1$ belongs to the pre-Hilbert space $\mathscr{C}_0/\mathcal{N}$ 
if and only if $n-m>2$.
\end{prop}

\begin{proof}
Recall that the reproducing kernel has a series expansion
$$\det(I_n-zw^*)^{-1}=\sum_{k=0}^{\infty} k!K_{k}(z,w),$$
where $K_{k}(z,w)$ is the reproducing kernel with respect to the Fock-Fischer norm for the $K$-type indexed by $k$. 
The functions 
\begin{eqnarray*}
z \mapsto \int_{Y}K_{2k}(z,v)d\omega
\end{eqnarray*}
are then $L$-invariant vectors in the $K$-type $2k$ and hence differ
from the $L$-invariants $\psi_k$ by some constants depending on $k$. We determine these by computing the integrals for a suitable choice of
$z$.

Before we begin with the computations, consider the fibration 
\begin{eqnarray*}
p: Y \rightarrow S^{n-1},
p(v)=v(e_1).
\end{eqnarray*} 
For $u \in S^{n-1}$, the fibre $p^{-1}(u)$ can be identified with the set of all rank $m-1$ partial isometries from
$\mathbb{R}^m$ to $(\mathbb{R}u)^{\perp}$. Moreover, $p$ is equivariant with respect to the actions of $O(n)$ on $Y$ and $S^{n-1}$.
Hence the equality
\begin{eqnarray}
\int_{S^{n-1}}fd\sigma=\int_Y f \circ p \,\,d\omega, \label{framskjutning}
\end{eqnarray} 
where $\sigma$ denotes the normalised rotation invariant measure on $S^{n-1}$, holds for all $f \in C(S^{n-1})$.

Choose now $z=\lambda e_1$, where $0 < \lambda <1$. Since $zv^t$ is a matrix of rank one, $\det(I_n-zv^t)^{-1}=(1-\mbox{tr}(zv^t))^{-1}$.
Hence 
\begin{eqnarray*}
\int_Y(1-\mbox{tr}(zv^t))^{-1}d\omega=\int_Y(1-\lambda v_{11})^{-1}d\omega=\int_Y(1-\lambda p(v)_1)^{-1}d\omega.
\end{eqnarray*} 
By \eqref{framskjutning}, we have
\begin{eqnarray*}
\int_Y(1-\lambda p(v)_1)^{-1}d\omega=\int_{S^{n-1}}(1-\lambda u_1)^{-1}d\sigma(u).
\end{eqnarray*}
Moreover, 
\begin{eqnarray*}
\int_{S^{n-1}}(1-\lambda u_1)^{-1}d\sigma(u)=\sum_{j=0}^{\infty}\lambda^j \int_{S^{n-1}}u_1^jd\sigma(u).
\end{eqnarray*}
The integrands on the right hand side depend only on the first coordinate, and hence the integrals can be written as integrals over 
the open interval $(-1,1)$ in $\mathbb{R}$ (cf. \cite{rudinft} 1.4.4.). In fact,
\begin{eqnarray*}
\int_{S^{n-1}}u_1^jd\sigma(u)=\frac{\Gamma(n/2)}{\Gamma(1/2)\Gamma((n-1)/2)} \int_{-1}^1(1-x^2)^{(n-2)/2-1}x^jdx.
\end{eqnarray*}
This integral is zero for odd $j$, and for $j=2k$, we have
$$\int_{-1}^1(1-x^2)^{(n-2)/2-1}x^jdx=B\left(\frac{2k+1}{2},\frac{n-1}{2}\right):=\frac{\Gamma(\frac{2k+1}{2})\Gamma(\frac{n-1}{2})}
{\Gamma(\frac{2k+n}{2})}.$$
Therefore, 
\begin{eqnarray*}
\int_{S^{n-1}}(1-\lambda u_1)^{-1}d\sigma(u)=\sum_{k=0}^{\infty}\frac{\Gamma(\frac{n}{2})\Gamma(\frac{2k+1}{2})}
{\Gamma(\frac{1}{2})\Gamma(\frac{2k+n}{2})} \lambda^{2k}. 
\end{eqnarray*}
From this, it follows that for an arbitrary  $z \in \mathscr{D}$, we have the expansion
\begin{eqnarray}
&&\int_{Y}\det(I_n-zv^t)^{-1}d\omega(v) \label{K-exp}\\
&&=\sum_{k=0}^{\infty}\frac{\Gamma(\frac{n}{2})\Gamma(\frac{2k+1}{2})} 
{\Gamma(\frac{1}{2})\Gamma(\frac{2k+n}{2})\Gamma(2k+1)} \psi_k(z). \nonumber
\end{eqnarray}
Since the functions $\psi_k$ are $L$-invariant, they are constant on the set $\{ru | u \in Y, 0<r<1\}$. This value equals
\begin{eqnarray}
\psi_k(ru)=r^{2k}\psi_k(u)=r^{2k}4^kk!\left(\frac{m}{2}\right)_k. \label{psi-ett}
\end{eqnarray}
Suppose now that $|r-r'|<\epsilon$. 
By \eqref{K-exp} and \eqref{psi-ett}, 
\begin{multline*}
\lefteqn{\int_{Y}K_1(ru,v)d\omega(v)-\int_{Y}K_1(r'u,v)d\omega(v)}\\
=\sum_{k=0}^{\infty}\frac{\Gamma(\frac{n}{2})\Gamma(\frac{2k+1}{2})} 
{\Gamma(\frac{1}{2})\Gamma(\frac{2k+n}{2})\Gamma(2k+1)} 4^kk!\left(\frac{m}{2}\right)_k(r^{2k}-(r')^{2k}),
\end{multline*}
and hence we have the estimate
\begin{multline}
\lefteqn{\left|\int_{Y}K_1(ru,v)d\omega(v)-\int_{Y}K_1(r'u,v)d\omega(v)\right| \leq}\\
\epsilon \sum_{k=0}^{\infty}\frac{\Gamma(\frac{n}{2})\Gamma(\frac{2k+1}{2})} 
{\Gamma(\frac{1}{2})\Gamma(\frac{2k+n}{2})\Gamma(2k+1)} 4^kk!\left(\frac{m}{2}\right)_k. \label{estimate}
\end{multline}
Applying Sterling's formula to the $k$th term on the right hand side yields
\begin{eqnarray}
\frac{\Gamma(\frac{n}{2})\Gamma(\frac{2k+1}{2})} 
{\Gamma(\frac{1}{2})\Gamma(\frac{2k+n}{2})\Gamma(2k+1)} 4^kk!\left(\frac{m}{2}\right)_k=O(k^{-\frac{n-m}{2}}).
\end{eqnarray}
Hence, the sum in \eqref{estimate} converges if and only if $n-m>2$. In this case, the corresponding
net $\left\{\int_Y K_1(r \cdot, v)d\omega\right\}_r$ is Cauchy in the supremum norm, and hence converges uniformly.
\end{proof}

\begin{lemma}\label{Ind-identifikation}
Consider the representation $\tau$ in \eqref{tauverkan}. On the space of continuous functions on $Y$, it is equivalent to the
representation $\mbox{Ind}_P^H(1 \otimes (i\lambda+\rho) \otimes 1)$, where 
$P$ is the minimal parabolic subgroup defined by the maximal abelian subspace $\mathfrak{a} \subset \mathfrak{p}$,
and $\lambda \in (\mathfrak{a}^{\mathbb{C}})^*$ is defined as
\begin{eqnarray}
-(i\lambda+\rho)|_{\mathfrak{a}_0}&=&-\frac{2(n-2)}{n-1}\rho_0, \label{killperp1}\\
-(i\lambda+\rho)|_{\mathfrak{a}_0^{\perp}}&=&0.\label{killperp2}
\end{eqnarray}
In fact, when the continuous functions on $Y$ are identified with right $L \cap M_0$-invariant functions on $L$, we can extend
them to functions on $H$ in such a way that the two representations are equal in this setting.
\end{lemma}

\begin{proof}
By \eqref{rho-jacobian}, we can rewrite the action of $H$ in \eqref{tauverkan} as
\begin{eqnarray}
\tau(h)f(x)=e^{-\frac{2(n-2)}{n-1}\rho_0(\log A_0(g^{-1}x)}f(\kappa(g^{-1}x)), 
\end{eqnarray}
where 
\begin{eqnarray*}
g^{-1}x=\kappa(g^{-1}x)m_0(g^{-1}x)A_0(g^{-1}x)n_0(g^{-1}x) \in LM_0A_0N_0.
\end{eqnarray*}
We now let $\lambda \in (\mathfrak{a}^\mathbb{C})^*$ be defined by the requirements \eqref{killperp1} and \eqref{killperp2}.

By \eqref{killperp2}, $-(i\lambda+\rho)$ has to annihilate all the restricted root spaces $\mathfrak{h}_{E_j^*-E_i^*}$, and hence
be of the form $c(E_1^*+\cdots+E_m^*)$ for some constant $c$. By $\eqref{rho}$ it follows that $c=-m(n-2)$.

Consider now the parabolically induced representation 
\begin{eqnarray*}
\mbox{Ind}_P^H(1 \otimes \exp(i\lambda+\rho) \otimes 1)
\end{eqnarray*}
acting on on continuous functions
on $H$.
By definition, this representation is defined on the space of continuous functions $f: H \rightarrow \mathbb{C}$ having
the $P$-equivariant property
\begin{eqnarray}
f(xman)=e^{-i(\lambda+ \rho)(\log a)}f(x). \label{P-ekvi}
\end{eqnarray}
The action of $H$ is given by
\begin{eqnarray}
f \stackrel{h}{\mapsto}e^{-(i\lambda+\rho)(A(h^{-1}x))}f(\kappa(h^{-1}x)). 
\end{eqnarray}

On the other hand, the restriction of the representation $\tau$ to the space of continuous functions on $Y$ coincides
with the $H$-action defined by the parabolically induced representation $\mbox{Ind}_{P_0}^H(\exp)$. Since $P \subset P_0$, and
\begin{eqnarray}
e^{-(i\lambda+\rho)(\log A(x))}=e^{-(i\lambda+\rho)(\log A_0(x))},
\end{eqnarray}
it follows that 
\begin{eqnarray}
\tau(h)f(x)=e^{-(i\lambda+\rho)(\log A(h^{-1}x))}f(\kappa(h^{-1}x), 
\end{eqnarray}
where $f$ is the extension of a continuous function on $Y$ to a $P_0$-equivariant function on $H$.
This finishes the proof.
\end{proof}

\begin{prop}
The operator $T: \mathscr{C}_0 \rightarrow \mathcal{O}(\mathscr{D})$ defined by
$$Tf(z)=\int_Y f(v)\det(I_n-zv^t)^{-1}d\omega(v)$$ 
is $H$-equivariant.
\end{prop}

\begin{proof}
We have
\begin{eqnarray*}
T(\tau(h)f)(z)&=&\int_Y J_h(h^{-1}v)^{\beta}f(h^{-1}v)K_1(z,v)d\omega\\
&=&\int_Y J_h(s)^{\beta+\frac{n-1}{n+m}}f(s)K_1(z,hs)d\omega\\
&=&J_h(h^{-1}z)^{-\frac{1}{n+m}}\int_Y J_h(s)^{\beta+\frac{n-1}{n+m}-\frac{1}{n+m}}f(s)K_1(h^{-1}z,s)d\omega\\
&=&\pi_1(h)(Tf)(z).
\end{eqnarray*}

\end{proof}

\begin{cor} \label{infkar}
The function $T1$ is a joint eigenfunction for all operators $\pi_1(Z)$, $Z \in Z(U(\mathfrak{h}^{\mathbb{C}}))$.
In particular, it is an eigenfunction for the Casimir operator, $\pi_1(\mathcal{C})$, with eigenvalue $-\frac{m(n-2)}{4}$. 
\end{cor}

\begin{proof}
By Lemma \ref{Ind-identifikation}, we can identify the extension of constant function 1 on $Y$ to a function on $H$ with the 
\emph{Harish-Chandra e-function} 
$e_{\lambda}:H \rightarrow \mathbb{C}$ given by
\begin{eqnarray}
e_{\lambda}(h)=e^{-i(\lambda+\rho)(\log A(h)}.
\end{eqnarray}
Moreover, the representation $\mbox{Ind}_P^H(1 \otimes \exp(i\lambda+\rho) \otimes 1)$ 
has infinitesimal character $i\lambda+\rho$ (cf. \cite{knapp1}, Ch. VIII). The value of the Casimir element
is $-(i\lambda+\rho)(\mathcal{C})=-(\langle \lambda,\lambda \rangle + \langle \rho, \rho \rangle)=-\frac{m(n-2)}{4}$ (cf. \cite{knapp-beyond}, Ch. V).
\end{proof}

\begin{prop}
The function $(T1)(z)=\int_{Y}\det(I_n-zv^t)^{-1}d\omega(v)$ belongs to $\mathscr{H}_1$.
\end{prop}
 
\begin{proof}
We rewrite the series expansion in \eqref{K-exp} using the orthonormal basis $\{\varphi_k\}$, i.e., 
\begin{eqnarray}
\int_{Y}K_1(z,v)d\omega(v)=\sum_{k=0}^{\infty} \alpha_k \varphi_k(z),
\end{eqnarray}
where $\alpha_k=\frac{\Gamma(\frac{n}{2})\Gamma(\frac{2k+1}{2})4^kk!\left(\frac{n}{2}\right)_{2k}^{1/2}\left( \frac{m}{2}\right)_{2k}^{1/2}} 
{\Gamma(\frac{1}{2})\Gamma(\frac{2k+n}{2})\Gamma(2k+1)}$.
By Sterling's formula
\begin{eqnarray}
 \alpha_k^2=O(k^{-(n-m)/2}),
\end{eqnarray}
and hence the series $\sum_k \alpha_k^2$ converges if and only if $n-m>2$.
\end{proof}
The operator $T$ maps the $H$-span (the set of all finite linear combinations 
$c_1 \tau(h_1)1+ \cdots + c_N \tau(h_1)1, h_i \in H, c_i \in \mathbb{C}$) of the function $1$ into $\mathscr{H}_1$. We introduce 
the temporary notation $H \cdot 1$ to denote this subspace. Moreover, we let $\mathcal{N}_{H \cdot 1}:=\mathcal{N} \cap H \cdot 1$.
\begin{prop}\label{isometrisk}
The equality
\begin{eqnarray}
\langle Tf,Tf \rangle_1=\langle f,f \rangle_{\mathscr{C}_0}
\end{eqnarray}
holds for $f \in H \cdot 1$.
\end{prop}

\begin{proof}
For $f \in H \cdot 1$ and $r \in (0,1)$, consider the function $Tf(r \cdot)$. 
We have
\begin{eqnarray}
Tf(rz)&=&\int_Y f(v)K_1(rz,v)d\omega(v)\\
&=&\int_Y f(v)K_1(z,rv)d\omega(v).
\end{eqnarray}
The square of the $\mathscr{H}_1$-norm is then given by
\begin{eqnarray*}
\|Tp(r \cdot)\|_1^2=\int_Y \int_Yf(\zeta)\overline{f(\eta)}K_1(r\zeta,r\eta)d\omega(\zeta)d\omega(\eta).
\end{eqnarray*}
These norms are uniformly bounded in $r$, and hence there is a convergent sequence $\{Tf(r_k \cdot)\}_k$ with respect to the
$\mathscr{H}_1$-norm. Since point evaluation functionals are continuous, we also have pointwise convergence, and hence this
limit function is $Tf$. Therefore,
\begin{eqnarray*}
\|Tf\|_1^2&=&\lim_{k \rightarrow \infty}\|Tf(r_k \cdot )\|_1^2\\
&=&\lim_{r \rightarrow 1}\int_Y \int_Y f(\zeta)\overline{f(\eta)}K_1(r^2\zeta,\eta)d\omega(\zeta)d\omega(\eta)\\
&=&\|f\|_{\mathscr{C}_0}^2.
\end{eqnarray*}
\end{proof}
We let $T_1$ denote the restriction of the operator $T$ to the subspace $H \cdot 1$. Then, we have the following corollary.
\begin{cor}
For the operator $T_1: H \cdot 1 \rightarrow \mathscr{H}_1$, 
\begin{eqnarray}
\ker\,T_1=\mathcal{N}_{H \cdot 1}.
\end{eqnarray}
\end{cor}
The operator $T_1$ then descends to an operator $U_1: H \cdot 1/\mathcal{N}_{H \cdot 1} \rightarrow \mathscr{H}_1$.
Now let $\mathcal{H}$ denote the Hilbert space completion of the space $H \cdot 1/\mathcal{N}_{H \cdot 1}$. We keep the letter
$\tau$ to denote the representation of $H$ of this space (in reality, the representation we mean is derived from $\tau$ by
first restricting, then descending to a quotient, and, finally, by extending uniquely to a Hilbert space completion).
\begin{prop}\label{irrkvot}
The representation $\tau$ of $H$ on $\mathcal{H}$ is irreducible.
\end{prop}

\begin{proof}
The representation $\tau$ is $H$-cyclic with a spherical ($L$-invariant) vector. Hence, there exists a unitary, $H$-equivariant
direct integral 
decomposition 
\begin{eqnarray}
S: \mathcal{H} \rightarrow \int_{\Lambda}\mathcal{H}_{\lambda}d\mu(\lambda),
\end{eqnarray}
where $\Lambda$ is a subset of the bounded spherical functions (or rather, the functionals on $\mathfrak{a}$ that
parametrise them), $\mu$ is some measure on $\Lambda$, and $\mathcal{H}_{\lambda}$ is the canonical spherical unitary
representation corresponding to the spherical function $\phi_{\lambda}$. 
For each $\lambda$, we let $v_{\lambda}$ denote the canonical spherical vector in $\mathcal{H}_{\lambda}$.

Suppose now that $\tau$ is not irreducible, i.e., the set $\Lambda$ is not a singleton set.
Then, we can choose two disjoint open subsets $\Omega_1, \Omega_2$ of $\Lambda$.
We define vectors $s_1$ and $s_2$ in the Hilbert space $\int_{\Lambda}H_{\lambda}d\mu$ by
\begin{eqnarray*}
s_1(\lambda)=\left\{ \begin{array}{c}
v_{\lambda}, \qquad \mbox{if} \,\lambda \in \Omega_1\\
0_{\lambda},\qquad \mbox{otherwise}\\ \end{array} \right.,\\
s_2(\lambda)=\left\{ \begin{array}{c}
v_{\lambda}, \qquad \mbox{if} \,\lambda \in \Omega_2\\
0_{\lambda}, \qquad \mbox{otherwise}\\ \end{array} \right..
\end{eqnarray*}
The vectors $S^{-1}s_1$ and $S^{-1}s_2$ are then linearly independent spherical vectors in $\mathcal{H}$. But, clearly, the only
spherical vectors in $\mathcal{H}$ are the (cosets modulo $\mathcal{N}_{H \cdot 1}$ of the) constant functions; 
a contradiction.
\end{proof}
We are now ready to state a subrepresentation theorem. The proof follows from Prop. \ref{isometrisk}, the above corollary, 
and Cor. \ref{infkar}.
\begin{thm} \label{embedding}
The operator $U_1$ can be extended to an isometric $H$-intertwining operator 
\begin{eqnarray}
U: \mathcal{H} \rightarrow \mathscr{H}_1.
\end{eqnarray}
Its image is isomorphic to the spherical unitary representation corresponding to the discrete point 
$\left\{i\left(\frac{1}{2}-\frac{n-m}{4}\right)\right\}$
in the spectral decomposition for the Casimir operator $\pi_1(\mathcal{C})$.
\end{thm}

\end{document}